\documentclass[12pt]{amsart}
\usepackage{amsfonts,mathrsfs,cite}
\newtheorem{thm}{Theorem}[section]
\newtheorem{theorem}[thm]{Theorem}
\newtheorem{corollary}[thm]{Corollary}
\newtheorem{lemma}[thm]{Lemma}
\newtheorem{proposition}[thm]{Proposition}

\theoremstyle{definition}

\theoremstyle{remark}
\newtheorem{remark}[thm]{Remark}

\numberwithin{equation}{section}

\renewcommand\Re{\mbox{Re\,}}
\renewcommand\Im{\mbox{Im\,}}
\newcommand\dom{\mbox{dom\,}}
\newcommand\diag{\mbox{diag\,}}
\newcommand\wt{\widetilde}
\newcommand\wh{\widehat}

\newcommand\bC{{\mathbb C}}
\newcommand\bN{{\mathbb N}}
\newcommand\bR{{\mathbb R}}
\newcommand\bZ{{\mathbb Z}}

\newcommand\ba{{\mathbf a}}
\newcommand\bb{{\mathbf b}}
\newcommand\bd{{\mathbf d}}
\newcommand\tba{{\widetilde{\mathbf a}}}
\newcommand\tbd{{\widetilde{\mathbf d}}}
\newcommand\tbb{{\widetilde{\mathbf b}}}
\newcommand\bc{{\mathbf c}}
\newcommand\bs{{\mathbf s}}
\newcommand\rt{{\mathrm t}}

\newcommand\mC{{\mathbf C}}
\newcommand\mS{{\mathbf S}}
\newcommand\cM{{\mathscr M}}

\newcommand\al{\alpha}
\newcommand\eps{\varepsilon}
\newcommand\la{\lambda}
\newcommand\si{\sigma}

\newcommand\gt{\mathfrak{t}}

\begin{document}

\title[Eigenvalue asymptotics]{Eigenvalue asymptotics for
    Sturm--Liouville operators with singular potentials${}^{\dag}$}%
\author[R.~O.~Hryniv, Ya.~V.~Mykytyuk]%
    {Rostyslav O.~Hryniv and Yaroslav V.~Mykytyuk}%
\address{Institute for Applied Problems of Mechanics and Mathematics,
3b~Naukova st., 79601 Lviv, Ukraine and
Lviv National University, 1 Universytetska st., 79602 Lviv, Ukraine}%
\email{rhryniv@iapmm.lviv.ua, yamykytyuk@yahoo.com}%

\thanks{${}^{\dag}$The work was partially supported by Ukrainian Foundation
for Basic Research DFFD under  grant No.~01.07/00172.
}%


\date{\today}%
%\dedicatory{}%
%\commby{}%

%%%%%%%%%%%%%%%%%%%%%%%%%%%%%%%%%%%%%%%%%%%%%%%%%%%%%%%%%%%%

\begin{abstract}
We derive eigenvalue asymptotics for Sturm--Liouville operators with
singular complex-valued potentials from the space $W^{\al-1}_{2}(0,1)$,
$\al\in[0,1]$, and Dirichlet or Neumann--Dirichlet boundary conditions.
We also give application of the obtained results to the inverse
spectral problem of recovering the potential from these two spectra.
\end{abstract}
\maketitle
%%%%%%%%%%%%%%%%%%%%%%%%%%%%%%%%%%%%%%%%%%%%%%%%%%%%%%%%%%%%

\section{Introduction}

In this paper we shall study eigenvalue asymptotics for
Sturm--Liouville operators on the interval $[0,1]$ with distributional
potentials. Namely, we assume that $q$ is a complex-valued distribution
from the Sobolev space~$W^{\al-1}_2(0,1)$, $\al\in[0,1]$, and consider
an operator $T$ that (formally) corresponds to the differential
expression
\begin{equation}\label{eq:1.1}
    l(f) := -f'' + q f
\end{equation}
and, say, Dirichlet boundary conditions. Explicitly, the operator~$T$
is defined by the regularisation method that was suggested
in~\cite{AEZ} for the particular potential $q(x)=1/x$ and was developed
by Savchuk and Shkalikov in~\cite{SS} for the class of distributional
potentials from $W^{-1}_2(0,1)$. (Incidentally, in this situation the
form-sum~\cite{AK,BN} and the generalized sum~\cite{KK} methods yield
the same operator). We observe that the considered class of singular
potentials include Dirac $\delta$-type and Coulomb $1/x$-type
interactions that are widely used in quantum mechanics and mathematical
physics; see also~\cite{KM} for other physical models leading to
potentials from negative Sobolev spaces.

The regularisation method consists in rewriting~\eqref{eq:1.1} as
\begin{equation}\label{eq:1.2}
    l(f) = l_\si(f) := -(f' -\si f)' - \si f',
\end{equation}
where $\si$ is any distributional primitive of~$q$. We fix one such
primitive $\si\in L_2(0,1)$ in what follows and call the expression
$f'-\si f=:f^{[1]}$ the \emph{quasi-derivative} of the function~$f$.
The natural $L_2$-domain of $l_\si$ is
\[                      % 1.3
    \dom l_\si = \{f\in W^1_1(0,1) \mid f^{[1]} \in W^1_1(0,1),\
            l_\si(f)\in L_2(0,1)\},
\]
and we observe that for $f\in\dom l_\si$ the derivative $f'=\si f+f^{[1]}$
belongs to $L_2(0,1)$ (but need not be continuous), so that $\dom l_\si
\subset W^1_2(0,1)$.

In the present paper, we shall only focus on Sturm--Liouville
operators $T_{\mathrm D}=T_{\si,{\mathrm D}}$ and $T_{\mathrm
N}=T_{\si,{\mathrm N}}$ that are generated by $l_\si$ and the
Dirichlet and the Neumann--Dirichlet boundary conditions
respectively, although other boundary conditions can also be
treated in a similar manner (see, e.g., \cite{KM} for periodic
and~\cite{SS1} for general regular boundary conditions). In other
words, $T_{\mathrm D}$ and $T_{\mathrm N}$ are the restrictions
of~$l_\si$ onto the domains
\begin{equation}\label{eq:1.4}
\begin{aligned}          % 1.4
    \dom T_{\mathrm D} &= \{f\in\dom l_\si \mid f(0)=f(1) = 0\},\\
    \dom T_{\mathrm N} &= \{f\in\dom l_\si \mid f^{[1]}(0)=f(1) = 0\}.
\end{aligned}
\end{equation}

It is known~\cite{SS} that the operators $T_{\mathrm D}$ and $T_{\mathrm
N}$ are closed, densely defined and have discrete spectra tending to
$+\infty$. We denote by $\la_{n}^2$ (resp., $\mu_n^2$) the eigenvalues of
$T_{\mathrm D}$ (resp.,~$T_{\mathrm N}$) counted with multiplicities and
arranged by increasing of the real---and then, if equal, imaginary---parts
of~$\la_n$ (resp., $\mu_n$). For definiteness, we shall always take $\la_n$
and $\mu_n$ from the set
\begin{equation}\label{eq:1.5}
    \Omega:= \{z\in\bC \mid  -\pi/2 < \arg z \le \pi/2\}\cup\{z=0\}.
\end{equation}
If $\al=0$, i.e., if $\si\in L_2$, then the numbers $\la_n$ and
$\mu_n$ obey the asymptotics~\cite{HMinv,Sav,SS,SS1}
\begin{equation}\label{eq:1.6}
    \la_n = \pi n + \wt\la_n, \qquad
    \mu_n = \pi (n-\tfrac12) + \wt\mu_n,
\end{equation}
where $(\wt\la_n)_{n\in\bN}$ and $(\wt\mu_n)_{n\in\bN}$ are some
$\ell_2$-sequences. It is reasonable to expect that if $\si$
becomes smoother, then the remainders $\wt\la_n$ and $\wt\mu_n$
decay faster; for instance, if $\al=1$, i.e., if $q\in L_2(0,1)$,
then the classical result~(see, e.g., \cite[Theorem~3.4.1]{Ma} or
\cite[Theorem~2.4]{PT}) states that $\wt\la_n,\wt\mu_n = {\mathrm
O}(n^{-1})$. Thus the problem arises to characterize the decay of
$\wt\la_n,\ \wt\mu_n$ depending on $\al\in[0,1]$.

Our interest in the above problem has stemmed from the inverse
spectral theory for Sturm--Liouville operators with singular
potentials. Namely, we proved in~\cite{HMinv2} that, as soon as
the numbers $\la_n^2$ and $\mu_n^2$ are real, strictly increase
with $n$, interlace, and obey~\eqref{eq:1.6} with
$\ell_2$-sequences $(\wt\la_n)$ and $(\wt\mu_n)$, then there
exists a unique real-valued $\si\in L_2(0,1)$ such that
$\{\la_n^2\}$ and $\{\mu_n^2\}$ are spectra of the
Sturm--Liouville operators $T_{\si,{\mathrm D}}$ and
$T_{\si,{\mathrm N}}$ respectively with distributional potential
$q=\si'\in W^{-1}_2(0,1)$. It is reasonable to believe that if
$\wt\la_n$ and $\wt\mu_n$ decay faster, then $\si$ will be
smoother. For example, the classical result of
Marchenko~\cite[Theorem~3.4.1]{Ma} claims that if, under the above
assumptions, we have in addition
\begin{equation}\label{eq:1.5a}
    \wt\la_n = \frac{A}n + \frac{\wt\la_n'}{n},\qquad
    \wt\mu_n = \frac{A}n + \frac{\wt\mu_n'}{n}
\end{equation}
with real $A$ and $\ell_2$-sequences $(\wt\la'_n)$ and
$(\wt\mu'_n)$, then the corresponding potential $q$ is
in~$L_2(0,1)$ (thus $\al=1$) and $A=\tfrac1{2\pi}\int_0^1q$, cf.
also~\cite[Theorem 3.3.1]{Le}. It would be desirable to
``interpolate'' between $\al=1$ and $\al=0$ and solve the inverse
spectral problem for all intermediate $\al\in(0,1)$. The essential
step towards such project is to derive eigenvalue asymptotics for
Sturm--Liouville operators with potentials in
$W^{\al-1}_2(0,1)$---i.e., to treat the direct spectral problem.
And indeed, based on the results obtained here, we completely
solve the inverse spectral problem for Sturm--Liouville operators
with potentials in the scale $W^{\al-1}_2(0,1)$, $\al\in[0,1]$, in
our paper~\cite{HMinv4}.

Another motivation for this work is the recent papers~\cite{KM} and
\cite{SS1}, where similar questions are considered. In particular,
Kappeler and M\"ohr in~\cite{KM} found eigenvalue asymptotics for the
Dirichlet and periodic Sturm--Liouville operators with complex-valued
potentials that are \emph{periodic} distributions from the space
$W^{\al-1}_2(0,1)$, $\al\in(0,1]$. The Dirichlet eigenvalues $\la_n^2$
were proved there to obey the asymptotics
\begin{equation}\label{eq:1.8}
    \la_n^2 = \pi^2 n^2 + \wh q(0) - \frac{\wh q(-2n) + \wh q(2n)}{2}
                + \nu_n,
\end{equation}
where $\wh q(n)$ is the $n$-th Fourier coefficient of~$q$ and the
sequence $(\nu_n)$ belongs to $\ell_2^{2\al-1-\eps}$ with $\eps>0$
arbitrary (the weighted $\ell_p^s$ spaces are defined at the end of
Introduction); see~\cite{KM} for more precise formulations. The authors
performed the Fourier transform to work in the weighted $\ell_2$ spaces
rather than in the Sobolev spaces and then derive the estimates for the
resolvent that yielded the detailed localisation of the spectrum.

Savchuk and Shkalikov~\cite{SS1} considered Sturm--Liouville operators
with complex-valued potentials $q$ that are distributional derivatives
of functions $u\in L_2(0,1)$. They generalised the notion of the
Birkhoff regular boundary conditions to this singular case and, for
Birkhoff regular boundary conditions, found eigenvalue and
eigenfunction asymptotics by means of the \emph{modified Pr\"ufer
substitution}. For the particular case of Dirichlet boundary conditions
and the function $u$ that is either
\begin{itemize}
\item[(i)] of bounded variation over $[0,1]$, or
\item[(ii)] Lipschitz continuous on~$[0,1]$ with exponent $\al\in(0,1)$, or
\item[(iii)] from $W^\al_2(0,1)$, $\al\in[0,\tfrac12)$,
\end{itemize}
the authors found several terms of asymptotic expansions for
eigenvalues and eigenfunctions. In particular, for $u\in
W^\al_2(0,1)$, $\al\in(0,\tfrac12)$, they gave the formula
\begin{equation}\label{eq:1.9}
    \la_n = \pi n - \int_0^1 u(t)\sin(2\pi nt)\,dt + s_n
\end{equation}
with $(s_n)\in \ell_p^{2\al}$ for all $p>1$; the remainders $s_n$
are further specified to be equal to
\begin{equation}\label{eq:1.9a}
\begin{aligned}
    s_n &:= - \frac1{2\pi n}\int_0^1\si^2(t)\cos(2\pi nt)\,dt
             + \frac1{2\pi n}\int_0^1 \si^2(t)\,dt \\
        &\quad - 2 \int_0^1 \int_0^t u(t)u(s)
            \cos(2\pi nt)\sin(2\pi ns)\,ds\,dt + \rho_n
\end{aligned}
\end{equation}
with $(\rho_n)\in \ell_1^{2\al}$.

 We also
mention the papers~\cite{An1,An2,CTS,CM,HMcL,KL,MAL,Zh}, where the
eigenvalue asymptotics were studied for several other types of singular
potentials.

In the present paper, having in mind further concrete applications
to the inverse spectral theory, we do not aim at the utmost
possible generality. Instead, we confine ourselves to the
Dirichlet and Neumann--Dirichlet boundary conditions and
potentials from the Sobolev space scale $W^{\al-1}_2(0,1)$,
$\al\in[0,1]$, though the derived formula for the characteristic
functions can be used to treat any generalised Birkhoff regular
boundary conditions and the developed methods allow further
application to potentials $q=\si'$ with $\si$ of bounded variation
or Lipschitz continuous as in~\cite{SS1}.

For the Dirichlet eigenvalues, we generalise the related results
of~\cite{KM} in several ways: a wider class of potentials is treated
(no periodicity is assumed and $\al=0$ is included), $\eps$ is removed
in the relation $(\nu_n)\in\ell_2^{2\al-1-\eps}$ for $\nu_n$ as
in~\eqref{eq:1.8} (see Remark~\ref{rem:6.2}), and more terms of
asymptotic expansion are given. The asymptotic formulae for $\la_n^2$
derived here are basically the same as in~\cite{SS1} (see
Remark~\ref{rem:6.1}); however, we allow the case $\al\ge\tfrac12$ and
simultaneously treat the Neumann--Dirichlet case. We give a special
representation of the remainders (required for the inverse analysis) as
sine Fourier coefficients of some functions from the $W_2^s$-scale.
This, however, does not yield an optimal result for the remainders
$s_n$ of~\eqref{eq:1.9} in terms of the $\ell_p^s$-scale---roughly
speaking, it only implies that $(s_n)\in \ell_2^{2\al}$ (cf.~the above
mentioned results of~\cite{SS1}). We extract from $\rho_n$
of~\eqref{eq:1.9a} one additional term falling into $\ell_1^{2\al}$ and
show that the modified remainders $\tilde\rho_n$ form an
$\ell_\infty^{\gamma}$-sequence with $\gamma=\min\{3\al,1+\al\}$;
observe that this result is incomparable to the inclusion
$(\rho_n)\in\ell_1^{2\al}$ proved in~\cite{SS1}. Also~\cite{SS1} gives
the eigenfunction asymptotics, which we do not study here (though the
derived formula for the Cauchy matrix allows such an analysis).

Our main result describes the asymptotics of $\la_n$ and $\mu_n$
in the following way. Here and hereafter, we denote by $s_n(f)$
and $c_n(f)$ the $n$-th sine and cosine Fourier coefficient of a
function~$f$, i.e.,
\[
    s_n(f) = \int_0^1 f(x)\sin (\pi n x)\,dx, \qquad
    c_n(g) = \int_0^1 f(x)\cos (\pi n x)\,dx, \qquad n\in
        \mathbb Z_+;
\]
$V$ and $R$ will stand for the operators in $L_2(0,1)$ given by
\[
    (Vf)(x)=(1-2x)f(x), \qquad
    (Rf)(x)=f(1-x) .
\]

\begin{theorem}\label{thm:main0} Assume that $q\in W^{\al-1}_2(0,1)$ for
some $\al\in[0,1]$ and fix an arbitrary distributional primitive
$\si\in W^\al_2(0,1)$ of~$q$. Then there exists a function $\wt\si \in
W_2^{2\al}(0,1)$ such that
\begin{equation}\label{eq:1.10a}
\begin{aligned}
    \la_n &= \pi n - s_{2n}(\si) - s_{2n}(\wt\si),\\
    \mu_n &= \pi (n-\tfrac12) + s_{2n-1}(\si) + s_{2n-1}(\wt\si).
\end{aligned}
\end{equation}
\end{theorem}

A more detailed version of this theorem reads as follows.

\begin{theorem}\label{thm:main}
Assume that $q\in W^{\al-1}_2(0,1)$ for some $\al\in[0,1]$ and fix
an arbitrary distributional primitive $\si\in W^\al_2(0,1)$
of~$q$. Then there exists a function $\omega\in
W^{\gamma}_2(0,1)$, $\gamma:= \min\{3\al, 1+\al\},$  such that
\begin{equation}\label{eq:1.10}
\begin{aligned}
    \la_n &= \pi n + s_{2n}(\si^-) - s_{2n}(\si) c_{2n}(V\si)
         +s_{2n}(\omega),\\
    \mu_n &= \pi (n-\tfrac12) + s_{2n-1}(\si^+) - s_{2n-1}(\si)
c_{2n-1}(V\si) +s_{2n-1}(\omega),
\end{aligned}
\end{equation}
where
\[
        \si^\pm(x):= \pm \si(x) \mp \int_0^x \si^2(t) \, dt + \int_0^{1-x}
        \si(x+t)\si(t) \, dt.
\]
\end{theorem}

Theorem~\ref{thm:main} and~\cite{SS1} imply the following

\begin{corollary}\label{cor:1.3}
Under the assumptions (and in the notations) of Theorem~\ref{thm:main}
\begin{equation}\label{eq:1.11}
    \la_n = \pi n + s_{2n}(\si^-) - s_{2n}(\si) c_{2n}(V\si) +
        {\mathrm O}(n^{-\gamma}),
        \quad n\to \infty.
\end{equation}
\end{corollary}

In rough terms, our main result states that if the potential $q$
belongs to the Sobolev space $W^{\al-1}_2(0,1)$, then the remainders
$\wt\la_n$ and $\wt\mu_n$ are sine Fourier coefficients of a function
from~$W^\al_2(0,1)$. This statement agrees with the earlier known
results for $\al=1$ and $\al=0$ mentioned above (see also~\cite{Zh} for
the case where $q$ is the derivative of a function of bounded
variation, i.e., where $q$ is a signed measure and \cite{An1,An2} for
Sturm--Liouville operators in impedance form with impedance functions
from $W^1_p(0,1)$). Notice that Theorem~\ref{thm:main0} implies the
related results of the paper~\cite{KM}, though does not imply the
results of~\cite{SS1} for the operator $T_\mathrm{D}$ and
$\al\in[0,\tfrac12)$. However, Corollary~\ref{cor:1.3} gives better
than in~\cite{SS1} estimates for the Dirichlet eigenvalues in the
uniform norm (i.e., in the $\ell_\infty^s$-scale). Moreover, our
approach is completely different from those of the works~\cite{KM} and
\cite{SS1} and requires only minor effort to get the next terms in
asymptotic eigenvalue expansions for small $\al$.

The organisation of the paper is as follows. In
Section~\ref{sec:fact} we derive an equivalent factorised form for
the differential expression $l_\si$ of~\eqref{eq:1.2} that is more
convenient for our purposes. We use this factorised form in the
next section to derive an integral representation for the
characteristic functions of the operators~$T_{\mathrm D}$ and
$T_{\mathrm N}$, and in Section~\ref{sec:smooth} we show that the
integrand in this representation possesses the desired smoothness.
The problem is thus reduced to finding asymptotics of zeros of
certain entire functions, which we establish in
Section~\ref{sec:zeros} and then prove the main results and some
corollaries in Section~\ref{sec:main}. Several applications to the
inverse spectral problem are given in Section~\ref{sec:appl}.
Finally, Appendix~\ref{sec:app.A} provides some necessary facts
from the interpolation theory.

Throughout the paper (in addition to the above-introduced notations
$s_n(f)$, $c_n(f)$ for sine and cosine Fourier coefficients of a
function $f$ and $V$ and $R$ for the multiplication operator by $1-2x$
and the reflection operator about $x=\tfrac12$ respectively) $W^s_p$
with $p\in[1,\infty)$ and $s\in\bR$ will be a shorthand notation for
the Sobolev space $W^s_p(0,1)$; we shall also abbreviate $L_p(0,1)$ to
$L_p$. The norm in $W^s_2$ is denoted by $|\cdot|_s$ (see
Appendix~\ref{sec:app.A}). The space $\ell_p^s$ consists of sequences
$x=(x_n)_{n\in\bN}$ with
\[
    |x|_{\ell_p^s}:=  \Bigl(\sum_{n\ge1} n^{ps} |x_n|^p\Bigr)^{1/p}< \infty
\]
and is a Banach space (a Hilbert space for $p=2$) under the
norm~$|\cdot|_{\ell_p^s}$. For $p=\infty$ the norm above should be
taken as
    $|x|_{\ell_\infty^s}:= \sup_{n\ge1} |x_n|n^s$.

%%%%%%%%%%%%%%%%%%%%%%%%%%%%%%%%%%%%%%%%%%%%%%%%%%%%%%%%%%%%

\section{Reduction of $l_\si$ to the factorised form}\label{sec:fact}

We recall that, for a given potential $q\in W^{\al-1}_2$ with
$\al\in[0,1]$, we have defined the Sturm--Liouville operator
$T_{\mathrm D}$ (resp., ~$T_{\mathrm N}$) corresponding to the
differential expression~\eqref{eq:1.1} and Dirichlet (resp.,
Neumann--Dirichlet) boundary conditions as $T_{\mathrm D} f = l_\si(f)$
for $f\in \dom T_{\mathrm D}$ (resp., as $T_{\mathrm N}f=l_\si(f)$ for
$f\in \dom T_{\mathrm N}$). Here $\si$ is a fixed distributional
primitive of $q$, $l_\si$ is given by
\[                      % 2.1
    l_\si(f) := -(f'-\si f)' - \si f'
        = -\Bigl(\frac{d}{dx}+\si\Bigr)\Bigl(\frac{d}{dx}-\si\Bigr)f
          - \si^2 f,
\]
and $\dom T_{\mathrm D}$ and $\dom T_{\mathrm N}$ are described
in~\eqref{eq:1.4}.

In this section, we shall derive a representation for the differential
expression $l_\si$ in a slightly different form. Roughly speaking, the
claim is that the last summand ($-\si^2f$) in the above formula can be
removed by changing $\si$ appropriately. More precisely, given $\tau\in
L_2$, we denote by $m_\tau$ the differential expression
\[
    m_\tau(f):= -\Bigl(\frac{d}{dx}+\tau\Bigr)\Bigl(\frac{d}{dx}-\tau\Bigr)f
\]
considered on the natural $L_2$-domain
\[
    \dom m_\tau = \{f\in W^1_1 \mid f'-\tau f\in W^1_1,
            \ m_\tau(f)\in L_2\}.
\]
Our aim here is to show that (under some not very restrictive
assumption) $l_\si$ coincides with $m_\tau$ for a suitable choice
of $\tau$. See also~\cite{KPST} for similar results on the whole
axis.

Denote by $\gt_{\mathrm N}=\gt_{\si,{\mathrm N}}$ the quadratic form of the
operator~$T_{\mathrm N}=T_{\si,{\mathrm N}}$. Integration by parts gives
that, for all $f\in \dom T_{\mathrm N}$,
\[
    \gt_{\mathrm N}[f]:= \int_0^1 l_\si (f)\overline{f}
            = \int_0^1 (|f'|^2 - \si f \overline{f'} - \si f' \overline{f}).
\]
It is not difficult to show (see, e.g., \cite{BS,HMper}) that the
quadratic form
    \(
        \int_0^1( \si f \overline{f'} - \si f'\overline{f})
    \)
is relatively bounded with respect to the quadratic form
    $\int_0^1 (|f'|^2 + |f|^2)$
with relative bound zero; hence Theorem~VI.1.33 of~\cite{Ka} implies
that $\gt_{\mathrm N}$ is sectorial and that its domain is the same as
for the unperturbed case $\si=0$, i.e., that
\[
    \dom \gt_{\mathrm N} = \{ f \in W^1_2 \mid f(1)=0\}.
\]
(As an aside we notice that, if we start from the quadratic form
$\gt_{\mathrm N}$ and denote by $S_{\mathrm N}$ the corresponding sectorial
operator, then $S_{\mathrm N} = T_{\mathrm N}$; thus the form-sum method
and the regularisation method yield the same operator. Another consequence
of this equality is that $\dom T_{\mathrm N} \subset \dom \gt_{\mathrm N}
\subset W^1_2$; we also note that the inclusion $\dom T_{\mathrm N} \subset
W^1_2$ follows from the relation $\dom l_\si\subset W^1_2$, which was
explained in Introduction. Finally, similar statements are also true for
$T_{\mathrm D}$.)

In the next proposition we assume that the quadratic form $\gt_{\mathrm N}$
is \emph{strictly accretive}, i.e., that $\Re\, \gt_{\mathrm N} [f]>0$ for
all nonzero $f\in \dom\gt_{\mathrm N}$. Since $\gt_{\mathrm N}$ is
sectorial, this situation can be achieved by adding to $q$ a suitable
positive constant and thus is not very restrictive.

\begin{proposition}\label{pro:2.1}
Assume that $\si\in W^{\al}_2$, $\al\in[0,1]$, and that the
quadratic form~$\gt_{\mathrm N}$ is strictly accretive. Then there
exists a function $\tau\in W^\al_2$ such that $\tau-\si\in
W^1_1\cap W^{2\al}_2$, $(\tau-\si)(0)=0$, and $l_\si = m_\tau$.
Moreover, the function
\[
\tilde\phi(x):=\tau(x) - \si(x) +\int_0^x \si^2(t) \, dt, \qquad
x\in [0,1],
\]
belongs to $W_2^{\gamma}$ with $\gamma =\min \{3\al, 1+\al \}$.
\end{proposition}

\begin{proof}
We shall take $\tau$ in the form $u'/u$, where $u$ is any function
satisfying the equation $l_\si(u)=0$ and not vanishing anywhere in the
interval $[0,1]$. After we have proved that such an $u$ exists and that
$\tau$ is of the required smoothness, verification of the equality
$l_\si = m_\tau$ becomes an easy algebraic exercise (see below).

Denote by $u$ a solution of the equation $l_\si(f)=0$ satisfying the
initial conditions $u(0)=1$ and $u^{[1]}(0)=0$. We recall that by
definition the equality $l_\si(f)=0$ is equivalent to the following
first-order system:
\begin{equation}\label{eq:2.5}
    \frac{d}{dx} \binom{f}{f^{[1]}} =
    \begin{pmatrix} \si & 1 \\ - \si^2 & -\si \end{pmatrix}
                 \binom{f}{f^{[1]}}.
\end{equation}
Since the entries of the $2\times2$ matrix in~\eqref{eq:2.5} are
summable, this system enjoys the standard existence and uniqueness
properties. In particular, the solution~$u$ with the stated
initial conditions exists and is unique; moreover, both $u$ and
$u^{[1]}$ belong to $W^1_1$ and, a posteriori, $u\in W^1_2$.

We claim that $u$ does not vanish on~$[0,1]$. Assume the contrary, i.e.,
let there exist $x\in(0,1]$ such that $u(x)=0$. Then integration by parts
gives
\[                      % 2.6
    0 = \int_0^x l_\si(u) \overline{u} =
        \int_0^x (|u'|^2 - \si u \overline{u'} - \si u' \overline u).
\]
We denote by $v$ the function from $W^1_2$ that coincides with $u$ on
$[0,x]$ and equals zero on~$[x,1]$ and observe that the above equation
implies $\gt_{\mathrm N}[v]=0$. Recall that the quadratic form
$\gt_{\mathrm N}$ is strictly accretive by assumption; henceforth we must
have $v=0$, which is impossible in view of the equality $v(0)=u(0)=1$. The
derived contradiction proves that $u(x)\ne0$ for every $x\in[0,1]$.

Put now $\phi:=u^{[1]}/u$ and $\tau:=\phi+\si= u'/u$. The
function~$\phi$ is in $W^1_1$ and satisfies the equation
$\phi'=-(\phi + \si)^2$ and the initial condition~$\phi(0)=0$ (so
that also $(\tau-\si)(0)=0$), and therefore
\begin{equation}\label{eq:2.7}
    \phi(x) = - \int_0^x \phi^2(t)\,dt -2 \int_0^x \phi(t)\si(t)\,dt
                -\int_0^x \si^2(t)\,dt.
\end{equation}
It follows from Lemma~\ref{lem:A.2al} that $\phi\in W^\al_2$ and
then repeated application of Lemma~\ref{lem:A.2al} shows that the
right-hand side of~\eqref{eq:2.7} belongs to $W^{2\al}_2$, so that
$\phi\in W^{2\al}_2$. Next, equality~\eqref{eq:2.7} implies that
\begin{equation*}
  \tilde\phi(x)= - \int_0^x \phi^2(t)\,dt -2 \int_0^x \phi(t)\si(t)\,dt,
\end{equation*}
and henceforth $\tilde\phi\in W_2^{\gamma}$ by
Lemma~\ref{lem:A.2al} as claimed.

Take now $f\in \dom l_\si$; then $f'-\tau f = f^{[1]} - \phi f \in W^1_1$
and, using the identity $\phi' = - (\phi+\si)^2$, we find that
\begin{align*}
    \Bigl(\frac{d}{dx}+\tau\Bigr)\Bigl(\frac{d}{dx}-\tau\Bigr)f
        &= (f^{[1]})' - \phi' f - \phi f' + \tau f^{[1]}
                                - \tau \phi f\\
        &= (f^{[1]})' + (\phi+\si)^2 f - \phi f' + (\phi + \si)(f'-\si f)
                                - (\phi+\si)\phi f\\
        &= (f^{[1]})' + \si f' = - l_\si(f).
\end{align*}
This shows that $l_\si \subset m_\tau$. The reverse inclusion is
established analogously, and, as a result, we get $l_\si = m_\tau$.
The proposition is proved.
\end{proof}

%%%%%%%%%%%%%%%%%%%%%%%%%%%%%%%%%%%%%%%%%%%%%%%%%%%%%%%%%%%%

\section{Integral representation of the
            characteristic functions}\label{sec:int}

Assume that $\si\in L_2$ is such that the quadratic form
$\gt_{\si,{\mathrm N}}$ is strictly accretive (see
Section~\ref{sec:fact} for definitions). Then by
Proposition~\ref{pro:2.1} there exists $\tau\in L_2$ such that
$\phi:=\tau-\si\in W^1_1$, $\phi(0)=0$, and $l_\si = m_\tau$.
Consider the differential equation $m_\tau u = \la^2 u$, i.e.,
\begin{equation}\label{eq:3.1}
    - \Bigl(\frac{d}{dx}+\tau\Bigr)\Bigl(\frac{d}{dx}-\tau\Bigr)u = \la^2u.
\end{equation}
It can be written as a first order system
\begin{equation}\label{eq:3.2}
    \frac{d}{dx} \binom{u_1}{u_2} =
    \begin{pmatrix} \tau & 1 \\ - \la^2 & -\tau \end{pmatrix}
                 \binom{u_1}{u_2}
\end{equation}
with $u_1\equiv u$ and $u_2 \equiv u'-\tau u=:u^{[1]}_\tau$. For any
$a,b\in\bC$, there exists a unique solution $(u_1,u_2)^T$ of~\eqref{eq:3.2}
subject to the initial conditions $u_1(0)=a$, $u_2(0)=b$,
whence~\eqref{eq:3.1} has a unique solution $v$ satisfying the initial
conditions $v(0)=a$ and $v^{[1]}_\tau(0)=b$.

Denote by $s(\cdot,\la)$ and $c(\cdot,\la)$ the solutions of
equation~\eqref{eq:3.1} obeying the initial conditions
\[                          % 3.3
    s(0,\la)            = c^{[1]}_\tau(0,\la)=0, \qquad
    s^{[1]}_\tau(0,\la) = c(0,\la)           =1.
\]
Observe that in view of the equality $\phi(0)=0$ we have
$u_\tau^{[1]}(0)=u^{[1]}(0) - \phi(0)u(0)= u^{[1]}(0)$ for any $u\in
\dom l_\si=\dom m_\tau$. Therefore the numbers $\pm\la_n$ coincide with
the zeros of the entire even function $s(1,\cdot)$, while $\pm\mu_n$
coincide with those of $c(1,\cdot)$. In both cases multiplicities are
taken into account (i.e., if some $\la^2$ is an eigenvalue
of~$T_{\mathrm D}$ of algebraic multiplicity $m$, then $\la$ is a zero
of~$s(1,\cdot)$ of order~$m$, and similarly for $T_{\mathrm N}$). The
functions $c(1,\cdot)$ and $s(1,\cdot)$ are called the
\emph{characteristic functions} of the operators $T_{\mathrm N}$ and
$T_{\mathrm D}$ respectively.

Our aim in this section is to show that the characteristic functions
$c(1,\cdot)$ and $s(1,\cdot)$ allow integral representations of a special
form (see~\eqref{eq:3.24} below). We do this by deriving first a special
integral representation for the Cauchy matrix of system~\eqref{eq:3.2}.

To begin with, we notice that the matrix-valued function
\[                          % 3.4
    U(x)= U(x,\la) :=  \begin{pmatrix}
        c(x,\la)            & s(x,\la)\\
        c^{[1]}_\tau(x,\la) & s^{[1]}_\tau(x,\la)
            \end{pmatrix}
\]
satisfies the initial condition $U(0)=I$ (with $I=\diag (1,1)$)
and solves the equation
\begin{equation}\label{eq:3.5}
        U' = (A + \tau J) U,
\end{equation}
where
\[                          % 3.6
    A=A_\la:=\begin{pmatrix} 0 & 1 \\ - \la^2 & 0 \end{pmatrix}, \qquad
    J:= \begin{pmatrix} 1 & 0 \\ 0 & -1 \end{pmatrix}.
\]
In other words, $U$ is the Cauchy matrix of system~\eqref{eq:3.2}. Since
$\tau\in L_2$, equation~\eqref{eq:3.5} with the initial condition $U(0)=I$
is uniquely soluble and the solution $U$ belongs to~$W^1_2$ entrywise.

We shall, however, need a more explicit formula for the Cauchy
matrix~$U$. The standard method of variation of parameters yields the
equivalent integral equation for~$U$ in the form
\begin{equation}\label{eq:3.7}
    U(x) = e^{xA} + \int_0 ^x e^{(x-t)A}\,\tau(t)JU(t)\,dt,
\end{equation}
where the exponent $e^{xA}$ can be explicitly calculated as
\[                          % 3.8
    e^{xA} = \begin{pmatrix} \cos\la x & \tfrac1\la\sin\la x \\
                                -\la\sin\la x & \cos\la x\end{pmatrix}.
\]
Integral equation~\eqref{eq:3.7} can be solved by the method of
successive approximations; namely, with
\begin{equation}\label{eq:3.9}
    U_0(x):=e^{xA}\quad\text{and}\quad
    U_{n+1}(x) = \int_0^x e^{(x-t)A}\,\tau(t)J U_n(t)\,dt
    \quad\text{for}\quad n\ge0,
\end{equation}
the solution formally equals
 \(                             % 3.10
    \sum_{n=0}^\infty U_n.
 \)
Our next aim is to show that this series converges in a suitable topology
and that the sum is indeed the Cauchy matrix.

To this end we endow the space
$\cM_2:=M(2,\bC)$
of all $2\times2$ matrices
with complex entries with the operator norm $|\cdot|$ of the
Euclidean~$\bC^2$ space and denote by $W^1_2(\cM_2)$ the Sobolev space of
$\cM_2$-valued functions on the interval~$[0,1]$.

\begin{lemma}\label{lem:3.1}
The series $\sum_{n=0}^\infty U_n$ with $U_n$ given by~\eqref{eq:3.9}
converges in $W^1_2(\cM_2)$ to the Cauchy matrix~$U$.
\end{lemma}

\begin{proof}
Bearing in mind the identity $Je^{xA}=e^{-xA}J$ and using recurrent
relations~\eqref{eq:3.9}, we derive the formula
\begin{equation}\label{eq:3.11}
    U_n(x) = \int_{\Pi_n(x)} e^{(x-2\xi_n(\rt))A}\,
             \tau(t_1) \cdots \tau (t_n) J^n\, dt_1\dots dt_n,
\end{equation}
in which we have put        % 3.12
\begin{gather*}
    \Pi_n(x) = \{\rt:= (t_1,\dots,t_n)\in \bR^n \mid
                0 \le t_n \le \dots \le t_1\le x\},\\
        \xi_n(\rt) = \sum_{l=1}^n (-1)^{l+1} t_l.
\end{gather*}
Observe that $0\le \xi_n(\rt)\le x$ for $\rt\in \Pi_n(x)$; thus, denoting
by~$C$ the maximum of $\left|e^{xA}\right|$ over the interval~$[-1,1]$, we
get the estimate
\[                          % 3.13
  |U_n(x)|
    \le C \int_{\Pi_n(x)}|\tau(t_1)|\cdots|\tau(t_n)|\,dt_1\dots dt_n
    = \frac{C}{n!} \Bigl(\,\int_0^x |\tau|\,\Bigr)^n.
\]
Differentiating recurrent relations~\eqref{eq:3.9}, we find that
\begin{equation}\label{eq:3.14}
    U'_n(x) = A U_n(x) + \tau(x) J U_{n-1}(x),
\end{equation}
and hence, with $C_1 := C (2|A|^2+3)^{1/2}$,
\[                          % 3.15
    \|U_n\|_{W^1_2(\cM_2)}:=
        \left(\,\int_0^1 \bigl(|U'_n|^2 + |U_n|^2\bigr)\right)^{1/2}
        \le \frac{C_1}{(n-1)!}  \left(\int_0^1 |\tau|^2\right)^{1/2}\,
                                \left(\int_0^1 |\tau|\,\right)^{(n-1)}.
\]
This estimate justifies convergence of the series
          $\sum_{k=0}^\infty U_k$
in the $W^1_2(\cM_2)$-topology to some $\cM_2$-valued function $V$ obeying
the initial condition $V(0)=I$. Bearing in mind~\eqref{eq:3.14} and
differentiating this series term-by-term, we see that $V$ satisfies
equation~\eqref{eq:3.5} and thus indeed equals the Cauchy matrix $U$.
\end{proof}

Our next aim is to get an integral representation for $U(1)$ of a special
form. Upon change of variables
 \(                          % 3.16
    s = \xi_n(\rt), \ y_l = t_{l+1}, \ l=1,2,\dots,n-1,
 \)
we recast the integral in~\eqref{eq:3.11} for $x=1$ as
\[                          % 3.17
    U_n(1) = \int_0^1 e^{(1-2s)A} \,\tau_n(s) J^n\,ds.
\]
Here $\tau_1 \equiv \tau(s)$ and, for all $n\in\bN$,
\begin{equation}\label{eq:3.18}
    \tau_{n+1}(s) = \int_{\Pi^*_n(s)} \tau(s+\xi_n({\mathrm y}))\,
         \tau(y_1)\,\cdots\,\tau(y_n)\, dy_1 \dots dy_n
\end{equation}
with
\begin{equation}\label{eq:3.19}
    \Pi^*_n(s) =  \{{\mathrm y}=(y_1,\dots,y_n)\in \bR^n \mid
        0 \le y_n \le y_{n-1} \le \dots \le y_1
            \le  s+ \xi_n({\mathrm y}) \le 1\}.
\end{equation}
Using the Cauchy--Schwarz inequality and Fubini's theorem, we find
that for every $n\in\bN$ the function $\tau_n$ belongs to $L_2$
and that
\begin{align*}                          % 3.20
     |\tau_n|^2_0 &= \int_0^1 |\tau_n(s)|^2\,ds\\
        &\le \frac1{(n-1)!} \int_0^1\,ds
            \int_{\Pi^*_{n-1}(s)} |\tau(s+\xi_{n-1}({\mathrm y}))\,
                \tau(y_1) \cdots \tau(y_{n-1})|^2\,
                    dy_1 \dots dy_{n-1}\\
        &= \frac1{(n-1)!}  \int_{\Pi_n(1)}
            |\tau(t_1)|^2\cdots|\tau(t_n)|^2 \,dt_1\dots dt_n
        = \frac1{(n-1)!n!} |\tau|^{2n}_0.
\end{align*}
It follows that the series
 \(                         % 3.21
    \sum_{n=1}^\infty (\pm1)^n \tau_n
 \)
converges in $L_2$ to some function $\tau^\pm$; putting
$K=\diag\{\tau^+,\tau^-\}$, we arrive at the desired representation for
$U(1)$:
\[                          % 3.22
    U(1) = e^{A} + \int_0^1 e^{(1-2s)A}K(s)\,ds.
\]
Spelling out the first row of this matrix equality, we get the following
result.

\begin{theorem}\label{thm:3.4}
Assume that $\si \in L_2$ is such that the quadratic form
$\gt_{\si,{\mathrm N}}$ is strictly accre\-tive and denote by $\tau\in L_2$
the function of Proposition~\ref{pro:2.1}, for which $l_\si = m_\tau$. Then
the characteristic functions $c(1,\cdot)$ and $s(1,\cdot)$ of the operators
$T_{\mathrm N}$ and $T_{\mathrm D}$ equal
\begin{equation}\label{eq:3.24}
\begin{aligned}
    c(1,\la) &= \cos \la + \int_0^1 \tau^+(s) \cos\la(1-2s) \,ds,\\
    s(1,\la) &=\frac{\sin\la}{\la} +  \int_0^1 \tau^-(s)
                \frac{\sin\la(1-2s)}{\la}\,ds,
\end{aligned}
\end{equation}
where the $L_2$-functions $\tau^\pm$ are defined by
\begin{equation}\label{eq:3.25}
    \tau^\pm = \sum_{n=1}^\infty (\pm1)^n \tau_n
\end{equation}
with $\tau_1 \equiv \tau$ and $\tau_{n+1}$ given by~\eqref{eq:3.18} for
all~$n\in\bN$.
\end{theorem}

In the case where $\si$ (and thus $\tau$) belongs to $W^\al_2$
with some positive $\al$ the functions $\tau_n$ are also smoother.
We shall establish this fact in the following section.

%%%%%%%%%%%%%%%%%%%%%%%%%%%%%%%%%%%%%%%%%%%%%%%%%%%%%%%%%%%%

\section{Smoothness of the functions $\tau_n$}\label{sec:smooth}

\def\mf{{\mathbf f}}

The derivation of the integral representations for the characteristic
functions $c(1,\cdot)$ and $s(1,\cdot)$ in Section~\ref{sec:int} only used
the fact that $\tau$ belongs to $L_2$. If, however, the potential $q$ is a
distribution from $W^{\al-1}_2$ with some $\al\in(0,1]$, then $\tau\in
W^\al_2$ by Proposition~\ref{pro:2.1}, and we can expect that the
functions~$\tau_n$ of~\eqref{eq:3.18} also have some additional smoothness.
The aim of this section is to make this statement precise.

Fix a natural $n\ge2$ and consider an $n$-linear mapping $I_n: (L_2)^n \to
L_2$ that acts according to the formula~(cf.~\eqref{eq:3.18})
\[                          % 4.1
    I_n(\mf)(s) =  \int_{\Pi^*_{n-1}(s)}
        f_1(s+\xi_{n-1}({\mathrm y}))f_2(y_1)\,\cdots\,f_n(y_{n-1})\,
        dy_1 \dots dy_{n-1},
\]
where $\mf=(f_1,\dots,f_n)\in (L_2)^n$, the set $\Pi^*_{n-1}(s)$ is defined
by~\eqref{eq:3.19}, and $\xi_{n-1}({\mathrm y}):=\sum_{l=1}^{n-1}
(-1)^{l+1}y_l$. In particular, we see that
\[                          % 4.2
    \tau_n = I_n(\tau,\dots,\tau).
\]
First we show that, indeed, $I_n$ maps $(L_2)^n$ into $L_2$.

\begin{lemma}\label{lem:4.1}
For any $\mf=(f_1,\dots,f_n)\in (L_2)^n$ the function $g:= I_n(\mf)$
belongs to $L_2$ and
\[                          % 4.3
    |g|_0 \le \frac{1}{\sqrt{(n-1)!}}\prod_{l=1}^n |f_l|_0.
\]
\end{lemma}

\begin{proof}
Repeating the arguments of Section~\ref{sec:int} used to prove that
$\tau_n\in L_2$, we find that
\[                           % 4.4
    |g(s)|^2  \le \frac1{(n-1)!} \int_{\Pi^*_{n-1}(s)}
        |f_1(s+\xi_{n-1}({\mathrm y}))|^2|f_2(y_1)|^2\,\cdots\,|f_n(y_{n-1})|^2\,
         dy_1 \dots dy_{n-1}
\]
and thus
\[                           % 4.5
    |\,g\,|^2_0 \le \frac1{(n-1)!} \int_{\Pi_n(1)}
                |f_1(t_1)|^2|f_2(t_2)|^2\dots |f_n(t_n)|^2 dt_1\dots dt_n
            \le \frac1{(n-1)!} \prod_{l=1}^n |f_l|^2_0
\]
as required.
\end{proof}

Now we give an equivalent formula for $I_n$. Let $g=I_n(\mf)$ and let $h$
be an arbitrary function in $L_2$; then the $L_2$-scalar product
$(g,h)_{L_2}$ of $g$ and $h$ can be recast as
\begin{equation}\label{eq:4.6}
    (g,h)_{L_2} = \int_{\Pi_n(1)} f_1(t_1)\cdots f_n(t_n)
        \overline{h}\bigl(\xi_n(\rt)\bigr)\,dt_1\dots dt_n.
\end{equation}
Since such scalar products with $h$ from a total set in $L_2$ completely
determine the function $g$, we can use this formula to define the action of
the mapping $I_n$.

The main result of this section is contained in
Theorems~\ref{thm:4.2} and \ref{thm:4.2a} that show how $I_n$ acts
between the spaces $W^\al_2$.

\begin{theorem}\label{thm:4.2}
For $\al\in[0,1]$, let $\mf\in (W_2^\al)^n$ and $g:= I_n(\mf)$.
\begin{itemize}
 \item [(a)] For all $n\ge 2$ the function $g$ belongs to $W_2^{2\al}$
 and there exists a positive $C$ independent of\/ $\mf$ and $n$ such that, with
 $r:=\max\{0,n-5\}$, we have
 \begin{equation*}
    |g|_{2\al} \le \frac{C}{\sqrt{r!}}
                        \prod_{l=1}^{n}|f_l|_\al.
 \end{equation*}
 \item [(b)] For all $n\ge 4$ the function $g$ belongs to $W_2^{3\al}$
 and there exists a positive $C$ independent of\/ $\mf$ and $n$ such that, with
 $r:=\max\{0,n-7\}$, we have
 \begin{equation*}
    |g|_{3\al} \le \frac{C}{\sqrt{r!}}
                        \prod_{l=1}^{n}|f_l|_\al.
 \end{equation*}
\end{itemize}
\end{theorem}

\begin{proof}
Since the family $\{W_2^\al\}$, $\al\in\bR$, constitutes a Hilbert scale,
by virtue of Theorem~\ref{thm:polylin-interp} and Lemma~\ref{lem:4.1} it
suffices to prove the theorem only for $\al=1$. Observe also that statement
(a) for $n=2$ holds in view of~Lemma~\ref{lem:A.I2}, so that we may assume
that $n\ge3$.

We shall divide the proof of the case $\al=1$ and $n\ge3$ into several
steps and shall throughout denote by $C_k$ positive constants independent
of $\mf$ and~$n$.

\textbf{Step 1.} Let $n \ge 3$, $\mf:=(f_1,\dots,f_n)\in (L_2)^n$
and $f_n\in W_2^1$. We shall show that in this case the function
$g$ belongs to $W_2^1$ and
\begin{equation}\label{eq:4.7a}
    |g|_1 \le \frac{C_1}{\sqrt{(n-3)!}}|f_n|_1\prod_{j=1}^{n-1}|f_j|_0
\end{equation}
for some positive constant $C_1$.

Let $\phi$ be an arbitrary test function (i.e., a $C^{\infty}$
function with support in $(0,1)$); then by the definition of the
distributional derivative we have
\[                          % 4.8
    (g',\phi) = - (g,\phi').
\]
Integration by parts gives
\begin{align*}              % 4.9
    - \int_0^{t_{n-1}} f_n(t_n) &\phi'\bigl(\xi_n(\rt)\bigr)\,dt_n
        = (-1)^n \int_0^{t_{n-1}} f_n(t_n)
                \phi'_{t_n}\bigl(\xi_n(\rt)\bigr)\,dt_n\\
    &=(-1)^{n+1}\int_0^{t_{n-1}} f'_n(t_n)\phi\bigl(\xi_n(\rt)\bigr)\,dt_n\\
    &\quad + (-1)^n f_n(t_{n-1})\phi\bigl((\xi_{n-2}(\tilde{\tilde \rt})\bigr)
        + (-1)^{n+1} f_n(0)\phi\bigl((\xi_{n-1}({\tilde \rt})\bigr),
\end{align*}
where we have put $\tilde \rt = (t_1,\dots,t_{n-1})$ and
 $\tilde{\tilde \rt}=(t_1,\dots,t_{n-2})$. It follows
now that
\begin{align*}              % 4.10
  -(g,\phi') &= (-1)^{n+1} \bigl(I_n(f_1,\dots,f_{n-1},f'_n),\phi\bigr)\\
   &\quad+ (-1)^{n+1} f_n(0)\bigl(I_{n-1}(f_1,\dots,f_{n-1}),\phi\bigr)\\
   &\quad+ (-1)^n \bigl(I_{n-2}(f_1,\dots,f_{n-3}, \tilde{f}_{n-2}),
                                                        \phi\bigr),
\end{align*}
where $\tilde{f}_{n-2}:= A_+(f_{n-2},f_{n-1},f_{n})$ and the
multilinear mapping $A_+$ is given by
\begin{equation}\label{eq:4.11}
    A_+(h_1,h_2,h_3)(x) := h_1(x)
            \int_0^{x} h_2(y) h_3(y)\,dy.
\end{equation}
We observe that $A_+$ acts boundedly from $(L_2)^3$ into $L_2$. In
fact, if $h_1,h_2,h_3\in L_2$, then $A_+(h_1,h_2,h_3)$ is a
product of an $L_2$-func\-ti\-on $h_1$ and a $W^1_1$-function
$\int_0^x h_2h_3$ and thus is in~$L_2$; moreover,
\begin{equation}\label{eq:4.11a}
    |A_+(h_1,h_2,h_3)|_0 \le |h_1|_0 \max_x \Bigl|\int_0^x h_2h_3\Bigr|
        \le |h_1|_0|h_2|_0|h_3|_0.
\end{equation}
Hence
\begin{equation}\label{eq:4.12}
\begin{aligned}
    g' = (-1)^{n} I_n(f_1,\dots,f_{n-1},f'_n)
         &+ (-1)^{n} f_n(0) I_{n-1}(f_1,\dots,f_{n-1})\\
         &+ (-1)^{n+1} I_{n-2}(f_1,\dots,f_{n-3},\tilde{f}_{n-2})
\end{aligned}
\end{equation}
in the sense of distributions. Since the right-hand side of the above
equation belongs to $L_2$ by Lemma~\ref{lem:4.1}, we conclude that $g\in
W^1_2$.

Recall that $W^1_2$ is continuously embedded into $C[0,1]$ and
thus there is $C_2>0$ such that $\max_x|f| \le C_2|f|_1$ for all
$f\in W^1_2$. Applying~Lemma~\ref{lem:4.1} to~\eqref{eq:4.12} and
using the inequalities $|f_n(0)| \le C_2 |f_n|_1$ and
 $|\tilde f_{n-2}|_0\le |f_{n-2}|_0 |f_{n-1}|_0 |f_{n}|_0$
(recall~\eqref{eq:4.11a}), we arrive at estimate~\eqref{eq:4.7a}.

\textbf{Step 2.} Let $n \ge 3$, $\mf:=(f_1,\dots,f_n)\in (L_2)^n$
and $f_1\in W_2^1$. We shall show that then the function $g$
belongs to $W_2^1$ and
\begin{equation}\label{eq:4.16a}
    |g|_1 \le \frac{C_1}{\sqrt{(n-3)!}}|f_1|_1\prod_{j=2}^{n}|f_j|_0
\end{equation}
with the same $C_1$ as in \eqref{eq:4.7a}.

A direct verification shows that, with $R$ being the reflection
operator about $x=\tfrac12$,
\[                          % 4.20
    I_n(f_1,\dots,f_n) = I_n(Rf_n,\dots,Rf_1)
\]
if $n$ is even and
\[                          % 4.21
    I_n(f_1,\dots,f_n) = R I_n(Rf_n,\dots,Rf_1)
\]
if $n$ is odd. Since $R$ is unitary in $W^\al_2$ for every
$\al\in[0,1]$ (for the cases $\al=0$ and $\al=2$ this is evident
and for intermediate values follows by interpolation, see
Theorem~\ref{thm:A.interp}), the inclusion $g\in W_2^1$ and
estimate~\eqref{eq:4.16a} follow from the results of Step~1.

\textbf{Step 3.} Let $n \ge 3$, $\mf:=(f_1,\dots,f_n)\in (L_2)^n$
and $f_1, f_n\in W_2^1$. Using relation~\eqref{eq:4.12} and
bounds~\eqref{eq:4.7a} and \eqref{eq:4.16a}, we easily conclude
that the function $g$ belongs to $W_2^2$ and that
\begin{equation}\label{eq:4.25}
        |g|_2 \le \frac{C_3}{\sqrt{r!}}|f_1|_1 |f_n|_1
                \prod_{l=2}^{n-1} |f_l|_0
\end{equation}
with some $C_3>0$ independent of $\mf$ and $r:=\max\{0,n-5\}$. The
only thing to be justified is that for $n=3$ the function
$I_1(\tilde f_1)\equiv\tilde f_1= A_+(f_1,f_2,f_3)$ belongs to
$W_2^1$ and
\[
    |\tilde f_1|_1\le C_4 |f_1|_1|f_2|_0|f_3|_1.
\]
This, however, easily follows from the formula
\[
    (\tilde f_1)'(x) = f_1'(x)\int_0^x f_2(t) f_3(t)\,dt
            + f_1(x)f_2(x)f_3(x)
\]
showing that $\tilde f_1'$ belongs to $L_2$ and providing the
suitable estimate of its $L_2$-norm.

Formula~\eqref{eq:4.25} combined with the remarks made at the
beginning of the proof and the obvious inequality $|h|_0\le|h|_1$
completes the proof of statement~(a).

\textbf{Step 4.} Let $n \ge 4$, $\mf:=(f_1,\dots,f_n)\in (L_2)^n$
and $f_2\in W_2^1$. We shall show that then the function $g$
belongs to $W_2^1$ and the following identity holds:
\begin{equation}\label{eq:4.25a}
\begin{aligned}
    g' &= - I_n(f_1,f'_2, f_3, \dots,f_n)
            - I_{n-2}(A_+(f_1,f_2,f_3), f_4, \dots, f_n)\\
        &\quad    + I_{n-2}(A_-(f_3,f_2,f_1), f_4, \dots, f_n)
            + I_{n-2}(f_1,A_+(f_4,f_2,f_3), f_5, \dots, f_n),
\end{aligned}
\end{equation}
where  $A_+$ is the mapping of~\eqref{eq:4.11} and $A_-$ is given
by
\[
    A_-(h_1,h_2,h_3)(x) := h_1(x)
            \int_x^{1} h_2(y) h_3(y)\,dy.
\]

Given an arbitrary test function $\phi$ and integrating by parts,
we get
\begin{align*}
    - \int_{t_3}^{t_1} f_2(t_2)\phi'(\xi_n(\rt))\,dt_2
        &= \int_{t_3}^{t_1}
            f_2(t_2)\phi'_{t_2}(\xi_n(\rt))\,dt_2\\
        &= - \int_{t_3}^{t_1} f'_2(t_2) \phi(\xi_n(\rt))\,dt_2\\
        &\quad    + f_2(t_1) \phi(\xi_{n-2}(\tilde\rt))
            - f_2(t_3) \phi(\xi_{n-2}(\tilde{\tilde\rt})),
\end{align*}
where ${\tilde \rt}=(t_3,t_4,\dots,t_n)$ and
 $\tilde{\tilde \rt}=(t_1,t_4,\dots,t_n)$.
Substituting this relation into the expression for $-(g,\phi')$
(cf.~\eqref{eq:4.6}), after simple calculations we get
\begin{multline*}
    -(g, \phi') = - (I_n(f_1,f'_2, f_3, \dots,f_n), \phi)
            - (I_{n-2}(A_+(f_1,f_2,f_3), f_4, \dots, f_n),\phi)\\
           + (I_{n-2}(A_-(f_3,f_2,f_1), f_4, \dots, f_n),\phi)
            + (I_{n-2}(f_1,A_+(f_4,f_2,f_3), f_5, \dots, f_n),\phi),
\end{multline*}
which yields~\eqref{eq:4.25a}.

\textbf{Step 5.} Assume that $n\ge 4$, $\mf:=(f_1,\dots,f_n)\in
(W_2^1)^n$ and $g:=I_n(\mf)$. Observe that the multilinear
transformations $A_+$ and $A_-$ map continuously $(W_2^1)^3$ into
$W_2^1$, which can be verified by direct calculation or using
Lemma~\ref{lem:A.2al}.

It follows from~\eqref{eq:4.25a} by Step~3 of the proof that
$g'\in W_2^2$, so that $g\in W_2^3$. The required norm estimate
also follows from~\eqref{eq:4.25} and the continuity properties of
$A_+$ and $A_-$. Thus statement (b) of the theorem is justified
for the case $\alpha =1$ and consequently, by the Interpolation
Theorem~\ref{thm:polylin-interp}, for all $\alpha\in [0,1]$. The
theorem is proved.
\end{proof}

For $n=3$ we have a slightly worse result.

\begin{theorem}\label{thm:4.2a}
Assume that $\al\in[0,1]$, $\gamma:=\min\{3\al,1+\al\}$,
$\mf:=(f_1,f_2,f_3)\in (W_2^\al)^3$, and $g:=I_3(\mf)$. Then the
function $g$ belongs to $W_2^\gamma$ and there is a constant $C$
independent of $\mf$ such that
\[
    |g|_\gamma \le C |f_1|_\al |f_2|_\al |f_3|_\al.
\]
\end{theorem}

\begin{proof}
Using the definition of $I_3$ and changing the variables via
$t_1=y_1-y_2$, $t_2=y_2$, we arrive at the representation
\begin{align*}
    g(s) &= \int_{\Pi_2^*(s)}
            f_1(s+y_1-y_2)f_2(y_1)f_3(y_2)\,dy_1\,dy_2\\
         &= \int_0^{s}dt_2  f_3(t_2)
             \int_0^{1-s} f_1(s+t_1)f_2(t_1+t_2)\,dt_1.
\end{align*}

We now put
\begin{align*}
   g_1(s):= \tilde I_3(\mf)(s)
            &:= \int_0^s \,dt_2  f_1(t_2)
                \int_0^{t_2} f_2(t_1) f_3(t_1)\,dt_1\\
            &\quad - \int_0^s \,dt_2  f_3(t_2)
                \int_{t_2}^{1} f_1(t_1) f_2(t_1)\,dt_1
\end{align*}
and
\[
    h(s) := J_3(\mf) := I_3(\mf)(s) + \tilde I_3(\mf)(s)
\]
and show that the function $h$ belongs to $W^{3\al}_2$ and that,
moreover,
\begin{equation}\label{eq:4.26}
    |h|_{3\al} \le C_1 |f_1|_\al|f_2|_\al|f_3|_\al
\end{equation}
for some $C_1>0$ independent of~$\mf$. Since $J_3$ is a
multilinear mapping, in view of Interpolation
Theorem~\ref{thm:polylin-interp} and Lemma~\ref{lem:4.1} it
suffices to treat only the case $\al=1$.

Assume therefore that $\mf\in (W_2^1)^3$. Direct calculations show that
\begin{align*}
    h'(s)=g'(s) + g'_1(s) &= f_1(s)\int_0^{s} f_2(t) f_3(t)\,dt
            - f_1(1) \int_0^s f_2(1-s+t) f_3(t)\,dt\\
           &\quad + \int_0^s dt_2 f_3(t_2)\int_0^{1-s}
            f_1'(s+t_1)f_2(t_1+t_2)\, dt_1.
\end{align*}
Integrating by parts in the last integral, we arrive at the
relation
\begin{align*}
    h'(s) = -I_3(f_1,f'_2,f_3),
\end{align*}
so that $h'\in W_2^2$ by Step~3 of the proof of
Theorem~\ref{thm:4.2}. Also
\[
    |h'|_{2} \le C_2 |f_1|_1|f'_2|_0|f_2|_1,
\]
which yields estimate~\eqref{eq:4.26} for $\al=1$. By
interpolation, the results hold also for all intermediate
$\al\in[0,1]$.

Consider now the function $g_1$. Put
\[
    J(\psi_1,\psi_2) := \int_0^s \psi_1(t)\psi_2(t)\,dt
\]
for $\psi_1,\psi_2\in L_2$; then direct calculations show that
\[
    g_1 = J(f_1,\psi_{23}) - J(f_3,R\psi_{12}),
\]
where $\psi_{23}:= J(f_2,f_3)$ and $\psi_{12}:= J(Rf_1,Rf_2)$.
Applying Lemma~\ref{lem:A.2al} twice we find that
$\psi_{12},\psi_{23}\in W_2^{2\al}$,
 $g_1\in W_2^{\al+\min\{2\al,1\}}= W_2^\gamma$ and that
\[
   |g_1|_\gamma \le C_3 |f_1|_\al |\psi_{23}|_{\min\{2\al,1\}}
            + C_3 |f_3|_\al |\psi_{12}|_{\min\{2\al,1\}}
            \le C_4 |f_1|_\al |f_2|_\al |f_3|_\al.
\]
The theorem is proved.
\end{proof}

\begin{remark}
The statement of the previous theorem cannot be improved in the
sense that the exponent $\gamma=\min\{3\al,1+\al\}$ cannot be made
larger. This follows from the fact that the results of
Lemma~\ref{lem:A.2al} are sharp. The same statement holds also for
$\tau_3= I_3(\tau,\tau,\tau)$.
\end{remark}

\begin{corollary}\label{cor:4.3}
Assume that $\al\in(0,1]$, $\gamma:=\min\{3\al,1+\al\}$, $\tau\in
W_2^\al$ and $\tau_n:=I_n(\tau,\dots ,\tau)$. Then for every
$n\ge3$ the function $\tau_n$ belongs to $W^{\gamma}_2$ and,
moreover,
\[
    |\tau_n|_{\gamma} \le \frac{C}{\sqrt{r!}}|\tau|^n_\al
\]
with some constant $C>0$ and $r:=\max\{0,n-7\}$.
In particular, the functions $\tau^\pm$ in Theorem~\ref{thm:3.4}
have the form
$$\tau^\pm = \pm\tau + \tau_2 + \phi^\pm$$
with some $W^{\gamma}_2$-functions $\phi^\pm$.
\end{corollary}

%%%%%%%%%%%%%%%%%%%%%%%%%%%%%%%%%%%%%%%%%%%%%%%%%%%%%%%%%%%%

\section{Asymptotics of zeros of some entire functions}\label{sec:zeros}

As we have seen in the previous sections, the eigenvalue asymptotics is
completely determined by the asymptotics of zeros for entire functions of a
special form. The main result of this section shows how this asymptotics
can be calculated.

Assume that $f$ is an arbitrary function from $W^\al_2$ and put
\begin{align*}                      % 5.1
    F_{\mathrm{c}}(\la)&:= \cos\la + \int_0^1 f(x) \cos[\la(1-2 x)]\,dx, \\
    F_{\mathrm{s}}(\la)&:= \frac{\sin\la}{\la} + \int_0^1 f(x)
            \frac{\sin[\la(1- 2x)]}{\la}\,dx.
\end{align*}
These are even entire functions of $\la$; we denote by
$\xi_{2n-1}$ and $\xi_{2n}$, $n\in\bN$, zeros of $F_{\mathrm{c}}$
and $F_{\mathrm{s}}$ respectively from the set $\Omega$
of~\eqref{eq:1.5}. We repeat every zero $\la\ne0$ according to its
multiplicity, and if $\la=0$ is a zero of $F_\mathrm{c}$ or
$F_\mathrm{s}$ of order $2m$, then we repeat it $m$ times among
$\xi_{2n-1}$ or $\xi_{2n}$ respectively. We shall order $\xi_k$ so
that $\Re \xi_{2n+1} > \Re \xi_{2n-1}$ or $\Re \xi_{2n+1} = \Re
\xi_{2n-1}$ and $\Im \xi_{2n+1} \ge \Im \xi_{2n-1}$ and similarly
for $\xi_{2n}$.

It is known (cf.~\cite[Ch.~1.3]{Ma} and~\cite{HMinv}) that for $f\in L_2$
the numbers $\xi_n$ have the form
\[                  % 5.1
 \xi_n = \frac{\pi n}2 + \tilde \xi_n
\]
for some $\ell_2$-sequence $(\tilde\xi_n)_{n\in\bN}$ (in particular, the
remainders $\tilde\xi_n$ are the Fourier coefficients of some
$L_2$-function). It is reasonable to expect that if $f$ is smoother (say,
from~$W^\al_2$), then $\tilde\xi_n$ decay faster. This is precisely what
the following theorem states.

Recall that for an $L_2$-function $g$ we have denoted by $s_n(g)$
and $c_n(g)$ its sine and cosine Fourier coefficients
respectively, and by $V:L_2\to L_2$ the operator of multiplication
by the function $(1-2x)$.

\begin{theorem}\label{thm:5.1}
Assume that $\al\in (0,1]$, $\gamma = \min \{3\al, 1+\al \} $,
$f\in W^\al_2$ and that the numbers $\tilde\xi_n$ are defined as
above. Then there exists a function $g\in W^{\gamma}_2$ such that
\[              % 5.3
    \tilde\xi_n = s_n(f) -s_n(f)c_n(Vf) +s_n(g), \qquad n\in\mathbb N.
\]
In particular, $(\tilde\xi_n)$ is a sequence of sine Fourier
coefficients of some function in~$W_2^\al$.
\end{theorem}

Before proceeding with the proof of the theorem, we introduce the
following spaces.

For any $g\in L_2$, we denote by $\bc(g)$ and $\bs(g)$ the
sequences $\bigl(c_n(g)\bigr){}_{n=0}^\infty$ and
$\bigl(s_n(g)\bigr){}_{n=0}^\infty$ of its cosine and sine Fourier
coefficients respectively and put
\[
   \mC_{\al} := \{\bc(g) \mid g\in W^\al_2 \}, \qquad
   \mS_{\al} := \{\bs(g) \mid g\in W^\al_2 \}, \qquad \al\in [0,2].
\]
The lineals $\mC_{\al}$ and $\mS_{\al}$ are algebraically embedded
into $\ell_2(\bZ_+)$ and become Banach spaces under the induced
norms
\[
    \|\bc(g)\|_{\mC_{\al}} := \|g\|_{W^\al_2}, \qquad
    \|\bs(g)\|_{\mS_{\al}} := \|g\|_{W^\al_2}.
\]
For any $\ba, \bb \in\ell_2(\bZ_+)$ we shall denote by $\ba\bb$
the entrywise product of $\ba$ and $\bb$, i.e., the element of
$\ell_2(\bZ_+)$ with entries $(\ba\bb)_n:=a_n b_n$.

To establish Theorem~\ref{thm:5.1}, we shall essentially rely on
the following three lemmata. The first of them
(Lemma~\ref{lem:5.2}) is proved in~Appendix~\ref{sec:app.A}, and
the other two are simple corollaries of well known facts and thus
their proofs are omitted.

\begin{lemma}\label{lem:5.2}
Suppose that $\alpha, \beta\in [0,1]$ and $\ba \in \mC_{\al}$,
$\bb \in \mS_{\al}$, $\tba\in\mC_{\beta}$,  $\tbb\in\mS_{\beta}$.
Then $\ba\tba \in \mC_{\alpha+\beta}$, $\bb\tbb \in
\mC_{\alpha+\beta}$, $\ba\tbb \in \mS_{\alpha+\beta}$; moreover,
there exists a positive constant $\rho>0$ such that
\begin{equation}\label{eq.51}
\begin{split}
    \|\ba\tba\|_{\mC_{\alpha+\beta}}
        &\le \rho\|\ba\|_{\mC_{\alpha}} \|\tba\|_{\mC_{\beta}},\\
    \|\bb\tbb\|_{\mC_{\alpha+\beta}}
        &\le \rho\|\bb\|_{\mS_{\alpha}} \|\tbb\|_{\mS_{\beta}},\\
    \|\ba\tbb\|_{\mS_{\alpha+\beta}}
        &\le \rho\|\ba\|_{\mC_{\alpha}}\|\tbb\|_{\mS_{\beta}}.
\end{split}
\end{equation}
\end{lemma}

\begin{lemma}\label{lem:5.3}
For every $\al\in[0,1]$, the operator $Vf(x)=(1-2x)f(x)$ acts boundedly
in~$W_2^\al$.
\end{lemma}

The claim follows directly from the Interpolation
Theorem~\ref{thm:A.interp}.

\begin{lemma}\label{lem:5.4}
Suppose that $\alpha\in [0,1/2)$. Then the Hilbert space
\[
 H_{\al}:=\{\ba=(a_n) \in \ell_2(\bZ_+)  \mid  \quad a_0=0, \quad
        \sum_{n=1}^{\infty} n^{2\al} |a_n|^2  < \infty  \}
\]
with norm
    $\|\ba\|_{H_{\al}}:=
        \left(\sum_{n=1}^{\infty} n^{2\al} |a_n|^2 \right)^{1/2}
    $
coincides with the space $\mS_{\alpha}$, and the norms
$\|\cdot\|_{H_{\al}}$ and $\|\cdot\|_{\mS_{\alpha}}$ are
equivalent.
\end{lemma}

This is a corollary of a well-known fact about Fourier transforms
of spaces $W_2^\al$, see, e.g.,~\cite{KM,LM,SS1}.

\begin{proof}[Proof of Theorem~\ref{thm:5.1}]
Using the obvious relations
\begin{align*}             % 5.7a
    \cos \xi_{2n-1} &= (-1)^n \sin \tilde\xi_{2n-1},\\
                    \sin \xi_{2n}  &= (-1)^n \sin \tilde\xi_{2n},\\
    \cos[\xi_{2n-1}(1-2x)]&=
            (-1)^n \sin [\tilde\xi_{2n-1}(1-2x)- (2n -1)\pi x)],\\
    \sin[\xi_{2n}(1-2x)]  &=
            (-1)^n \sin [\tilde\xi_{2n}(1-2x)- 2n\pi x)]
\end{align*}
in the equalities $F_{\mathrm{c}}(\xi_{2n-1})=0$ and
$F_{\mathrm{s}}(\xi_{2n})=0$, we find that
\begin{equation}\label{eq:5.7}
    \sin \tilde{\xi}_n + \int_0^1 f(x) \sin \bigl[\tilde{\xi}_n(1-2x)-\pi nx \bigr]\,dx
          =0, \qquad n\in\bN.
\end{equation}

Writing $\sin \bigl[ \tilde{\xi}_n(1-2x)-\pi nx\bigr]$ as
 $\sin[\tilde{\xi}_n(1-2x)]\cos(\pi nx)- \sin(\pi nx)\cos[\tilde{\xi}_n(1-2x)]$, developing
 $\sin[\tilde{\xi}_n(1-2x)]$ and  $\cos[\tilde{\xi}_n(1-2x)]$
into the Taylor series, and then changing summation and
integration order (which is allowed in view of the absolute
convergence of the Taylor series and the integrals), we
represent~\eqref{eq:5.7} as
\begin{equation}\label{eq:5.8}
            \sin\tilde{\xi}_n   +
        \sum_{k=0}^\infty (-1)^k \frac{\tilde{\xi}_n^{2k+1}}{(2k+1)!}
                c_n(V^{2k+1}f)
    -   \sum_{k=0}^\infty (-1)^k \frac{\tilde{\xi}_n^{2k}}{(2k)!}
                s_n(V^{2k}f) =0.
\end{equation}

Set $\tilde\xi_0:=0$; then, as was mentioned above, the sequence
$\ba:=(\tilde\xi_n)_{n\in\mathbb Z_+}$ belongs to $\ell_2(\bZ_+)$, so that
$\ba\in \mS_0$. Define the sequence $\bd:=(d_n)_{n\in\bZ_+}$ through the
relation
\[
  \bd:=  \sum_{k=0}^\infty  \frac{(-1)^k}{(2k)!} \ba^{2k}\bs(V^{2k}f)
 -\sum_{k=0}^\infty  \frac{(-1)^k}{(2k+1)!} \ba^{2k+1} \bc(V^{2k+1}f).
\]
Using Lemmata~\ref{lem:5.2} and~\ref{lem:5.3} and denoting by
$\rho_1$ the norm of the operator $V$ in $W_2^\al$, we find that
\begin{align*}
    \|\ba^{2k}\bs(V^{2k}f)\|_{\mS_{\al}}
        &\le (\rho\rho_1)^{2k} \|\ba\|^{2k}_{\mS_0}\|f\|_{W_2^{\al}},\\
    \|\ba^{2k+1}\bc(V^{2k+1}f)\|_{\mS_{\al}}
        &\le(\rho\rho_1)^{2k+1}\|\ba\|^{2k+1}_{\mS_0}\|f\|_{W_2^{\al}},
\end{align*}
whence $\bd\in \mS_{\al}$.

Equation~\eqref{eq:5.8} implies that $\sin\tilde{\xi}_n =d_n$ for
all $n\in\bZ_+$; henceforth there exists $n_0\in\bN$ such that
\begin{equation}\label{eq:5.10}
  |\tilde{\xi}_n|\le 2|d_n|
\end{equation}
for all $n\ge n_0$. Fix an arbitrary number
    $\beta\in (\al/3,\al/2)$;
then $\beta < \al/2 \le 1/2$ and $\bd\in \mS_{\beta}$.
Lemma~\ref{lem:5.4} and inequalities~\eqref{eq:5.10} now yield the
inclusion $\ba\in \mS_{\beta}$.

Since the sequence $(\sin\tilde\xi_n)_{n\in\bZ_+}$ can be written
as
\[
    \sin\ba:=\sum_{k=0}^\infty (-1)^k \frac{\ba^{2k+1}}{(2k+1)!}
\]
and since by Lemma~\ref{lem:5.2} and the inequality $2\beta<1$ the
series
\[
    \tbd := \sin \ba - \ba =\sum_{k=1}^\infty (-1)^k \frac{\ba^{2k+1}}{(2k+1)!}
\]
converges in $\mS_{3\beta}$, we find that
\begin{equation}\label{eq:5.12}
 \ba=\bd -\tbd \in \mS_{\al} + \mS_{3\beta}\subset \mS_{\al}.
\end{equation}

By the definitions of $\bd$ and $\tbd$ equality~\eqref{eq:5.12}
can be recast as
\[
   \ba= \sum_{k=0}^\infty (-1)^k \frac{\ba^{2k}}{(2k)!} \bs(V^{2k}f)
    -\sum_{k=0}^\infty (-1)^k \frac{\ba^{2k+1}}{(2k+1)!} \bc (V^{2k+1}f)
    -\sum_{k=1}^\infty (-1)^k \frac{\ba^{2k+1}}{(2k+1)!}.
\]
Since in view of Lemmata~\ref{lem:5.2} and ~\ref{lem:5.3} the sum
 $$
\sum_{k=1}^\infty (-1)^k \frac{\ba^{2k}}{(2k)!} \bs(V^{2k}f)
 -\sum_{k=1}^\infty (-1)^k \frac{\ba^{2k+1}}{(2k+1)!} \bc (V^{2k+1}f)
 -\sum_{k=1}^\infty (-1)^k \frac{\ba^{2k+1}}{(2k+1)!}
 $$
falls into $\mS_{\gamma}$, we conclude that
\[
\ba - \bs(f) + \ba \bc (Vf) \in \mS_{\gamma},
\]
so that also
\[
   [\ba - \bs(f)] \bc (Vf) \in
         [-\ba \bc (Vf) + \mS_{\gamma} ]\bc (Vf) \subset \mS_{\gamma}.
\]
This yields the desired result
\[
  \ba - \bs(f) + \bs(f)\bc (Vf)\in \mS_{\gamma},
\]
and the theorem is proved.
\end{proof}

%%%%%%%%%%%%%%%%%%%%%%%%%%%%%%%%%%%%%%%%%%%%%%%%%%%%%%%%%%%

\section{Proof of Theorems~\ref{thm:main0},~\ref{thm:main}
and Corollary ~\ref{cor:1.3}}\label{sec:main}

It suffices to establish only the refined
asymptotics~\eqref{eq:1.10}. Indeed, since by
Lemmata~\ref{lem:A.2al} and \ref{lem:A.I2} the functions
 \(
    \int_0^x \si^2(t) \, dt
 \)
and
 \(
         \int_0^{1-x} \si(x+t)\si(t) \, dt
 \)
belong to $W_2^{2\al}$ and by Lemmata~\ref{lem:5.2} and
~\ref{lem:5.3} the numbers $s_{n}(\si)c_n(V\si)$ are $n$-th sine
Fourier coefficients of some function from $W_2^{2\al}$,
formula~\eqref{eq:1.10a} follows from~\eqref{eq:1.10}.

Assume first that the potential $q\in W^{\al-1}_2$ is such that the
associated quadratic form~$\gt_{\mathrm N}$ is strictly accretive. Then by
Proposition~\ref{pro:2.1} we can find a function $\tau\in W^\al_2$ such
that $\phi:= \tau-\si\in W^{2\al}_2\cap W^1_1$ and  $l_\si =  m_\tau$. By
Corollary~\ref{cor:4.3} the functions $\tau^\pm$ of~\eqref{eq:3.25} can be
represented as
\[
    \tau^\pm = \pm\tau + \tau_2 + \phi^\pm
\]
with $\tau_2=I_2(\tau,\tau)$ and some $W^{\gamma}_2$-functions
$\phi^\pm$. Taking into account Proposition~\ref{pro:2.1} and
Lemma~\ref{lem:A.I2}, we conclude that
\begin{equation}\label{eq:6.1}
    \tau^\pm = \si^\pm + \tilde\phi^\pm,
\end{equation}
where
\begin{equation}\label{eq:6.2}
    \si^\pm(x):= \pm \si(x) \mp \int_0^x \si^2(t) \, dt +
         I_2(\si,\si)(x)
\end{equation}
and $\tilde\phi^\pm \in W^{\gamma}_2$. By virtue of
Theorems~\ref{thm:3.4} and ~\ref{thm:5.1} there exist functions
$\psi^\pm\in W^{\gamma}_2$ such that
\begin{equation}\label{eq:6.3}
\begin{split}
    \wt \la_n    =   & s_{2n}(\tau^-) - s_{2n}(\tau^-) c_{2n}(V\tau^-)
         +s_{2n}(\psi^-), \\
    \wt \mu_n    =   & s_{2n-1}(\tau^+) - s_{2n-1}(\tau^+) c_{2n-1}(V\tau^+)
         +s_{2n-1}(\psi^+).
\end{split}
\end{equation}
Since  $(\tau^\pm \mp\si)\in  W_2^{2\al}$,
equalities~\eqref{eq:6.1}--\eqref{eq:6.3} and Lemma~\ref{lem:5.2}
yield the representation
\begin{equation}\label{eq:6.4}
\begin{split}
    \wt \la_n    =   & s_{2n}(\si^-) - s_{2n}(\si) c_{2n}(V\si)
         +s_{2n}(\wt\psi^-), \\
    \wt \mu_n    =  & s_{2n-1}(\si^+) - s_{2n-1}(\si) c_{2n-1}(V\si)
         +s_{2n-1}(\wt\psi^+)
\end{split}
\end{equation}
with some $\wt\psi^\pm\in W^{\gamma}_2$. In remains to put
$\omega:= \tfrac12(\wt\psi^+ + R\wt\psi^+ + \wt\psi^-
-R\wt\psi^-)$ to get the required formula.

In a generic situation we add a suitable constant~$C$ to the
potential~$q$ to get a potential $\wh q:= q+C$ that falls into the
above-considered case. (This can be done since the quadratic form
$\gt_{\mathrm N}$ is bounded below, see Section~\ref{sec:fact}.)
Then the corresponding Dirichlet eigenvalues $\wh \la^2_n$ and
Neumann--Dirichlet eigenvalues $\wh \mu^2_n$ have the
form~\eqref{eq:6.4} with $\si$ replaced by $\wh\si:=\si+Ct$ and
with $\wh\si^\pm$ calculated as in~\eqref{eq:6.2} for $\wh\si$
instead of~$\si$. Since $\si-\wh\si \in W_2^2$, it is easily seen
that $\si^\pm-\wh\si^\pm \in W_2^{\gamma}$. Calculating now the
integrals, we arrive at the relations
\begin{align*}
    \wh \la_n    &=    \pi n +s_{2n}(\si^-) - s_{2n}(\si) c_{2n}(V\si)
         +s_{2n}(\wh\psi^-), \\
    \wh \mu_n    &=  \pi(n-1/2) + s_{2n-1}(\si^+)
            - s_{2n-1}(\si) c_{2n-1}(V\si)
            + s_{2n-1}(\wh\psi^+)
\end{align*}
for some $\wh\psi^\pm\in W_2^\gamma$. Since $ \wh \la_n = \pi n +
a_n$ with $(a_n)\in\ell_2$, we find that
\[              % 6.3
  \wh\la_n -\la_n =\wh\la_n -\sqrt{\wh\la^2_n -C} =\frac{C}{2\wh\la_n} +
  {\mathrm O}(\wh\la_n^{-3})= \frac{C}{2\pi n} + \frac{b_n}{n^2}
\]
for some $(b_n)\in\ell_2$, so that there exists a function
$\chi\in W_2^2$ such that
\[
    s_{2n}(\chi)= \wh\la_n -\la_n, \qquad n\in \bN.
\]
Thus
\begin{align*}              % 6.3
  \wt\la_n &= s_{2n}(\si^-) - s_{2n}(\si) c_{2n}(V\si)
         +s_{2n}(\wh\phi^-),
\end{align*}
with some $\wh\phi^-\in W_2^{\gamma}.$ Similar arguments work for
$\wt\mu_n$ and yield the representation
\begin{align*}              % 6.3
  \wt\mu_n &= s_{2n-1}(\si^+) - s_{2n-1}(\si) c_{2n-1}(V\si)
         +s_{2n-1}(\wh\phi^+),
\end{align*}
with some $\wh\phi^+\in W_2^{\gamma}$. It remains to put
$\omega:=\tfrac12 (\wh\phi^+ + R\wh\phi^+ + \wh\phi^-
-R\wh\phi^-)$, and the proof of Theorem~\ref{thm:main} is
complete.

\begin{remark}\label{rem:6.1}
Straightforward calculations show that the asymptotic formula
\eqref{eq:1.10} established here corresponds to the asymptotic
formulae~\eqref{eq:1.9}--\eqref{eq:1.9a} of the paper~\cite{SS1} with
$\rho_n= 2s_{2n}(\si)c_{2n}(t\si) + s_{2n}(\omega)$. Since
$(\rho_n)\in\ell_2^{2\al}$ by the results of~\cite{SS1} and
$s_{2n}(\si)c_{2n}(t\si)$ falls into $\ell_2^{2\al}$ for
$\al\in[0,\tfrac12)$, one concludes that
$(s_{2n}(\omega))\in\ell_2^{2\al}$ for such $\al$.
Corollary~\ref{cor:1.3} states that, moreover, $(s_{2n}(\omega))\in
\ell_\infty^{\gamma}$ for all $\al\in[0,1]$.
\end{remark}

\begin{proof}[Proof of Corollary~\ref{cor:1.3}]
Assume that $\si\in W^\al_2(0,1)$ for some $\al\in[0,1]$ and
$\omega\in W^{\gamma}_2$ is the function of
Theorem~\ref{thm:main}. If $\al\in[0,1/3]$, then by
Lemma~\ref{lem:A.5} one has
\[
   s_{2n}(\omega)=\mathrm{O}(n^{-3\alpha}),  \qquad n\to\infty,
\]
and, in view of~\eqref{eq:1.10},
\[
    \la_n = \pi n + s_{2n}(\si^-) - s_{2n}(\si) c_{2n}(V\si)
       + \mathrm{O}(n^{-3\alpha}), \qquad n\to\infty.
\]

Let now $\al\in (1/3,1]$. Since the function $\si\in W_2^{\al}$ satisfy
then condition~(vi) of Theorem~3.12 of~\cite{SS1}, one has
\[
    \la_n = \pi n - s_{2n}(\si) + a_n
\]
with $(a_n)\in\ell_1$. In view of~\eqref{eq:1.10} this yields
\[
       s_{2n}(\omega)= - s_{2n}(\psi) - b_n,
\]
where
\begin{equation}\label{eq:6.6}
    \psi(x) := (\si^-+\si)(x)
     = \int_0^x \si^2(t) \, dt + \int_0^{1-x} \si(x+t)\si(t) \, dt
\end{equation}
and $(b_n)\in\ell_1$. By virtue of Lemmata~\ref{lem:A.2al}
and~\ref{lem:A.I2} the function $\psi$ belongs to the
space~$W_2^{2\alpha}$. Since, moreover, $\psi(0)=\psi(1)$,
Lemma~\ref{lem:A.5} (see also Remark~\ref{rem:A.sin-coef}) yields
the inclusion $(s_{2n}(\psi))_{n\in\mathbb N}
\in\ell_2^{2\al}\subset\ell_1$, i.e.,
$(s_{2n}(\omega))_{n\in\mathbb N} \in\ell_1$. Taking into account
that $\gamma\ge1$, we conclude that $\omega(0)=\omega(1)$.
Lemma~\ref{lem:A.5} yields now the inclusion
$(s_{2n}(\omega))\in\ell_2^\gamma$ and the required
relation~\eqref{eq:1.11} follows.
\end{proof}

\begin{remark}\label{rem:6.2}
Developing the above arguments, we can show that the function $\wt\si$ of
Theorem~\ref{thm:main0} is such that $(s_{2n}(\wt\si))\in\ell_2^{2\al}$.
For $\al\in[0,\tfrac12)$ this is shown in~\cite{SS1}. Observe that
\[
    s_{2n}(\wt\si)= s_{2n}(\omega) + s_{2n}(\si)c_{2n}(V\si)+s_{2n}(\psi)
\]
with the function $\omega$ of Theorem~\ref{thm:main} and $\psi$
of~\eqref{eq:6.6}. For $\al\in[\tfrac12,1]$ the proof of
Corollary~\ref{cor:1.3} establishes the inclusions
$(s_{2n}(\omega))\in\ell_2^\gamma\subset\ell_2^{2\al}$ and
$(s_{2n}(\psi))\in\ell_2^{2\al}$, while
\[
    \bigl(s_{2n}(\si)c_{2n}(V\si)\bigr) \in
    \ell_\infty^\al\cdot\ell_2^\al\subset\ell_2^{2\al}
\]
by Lemma~\ref{lem:A.5}.

If $q\in W^{\al-1}_2$ is periodic, then we find that
\[
    s_{2n}(\si) = -\frac{c_0(q)}{2\pi n}
            + \frac{c_{2n}(q)}{2\pi n}.
\]
Inserting this into~\eqref{eq:1.10a} and squaring, we get
\[
    \la_n^2 = \pi^2 n^2 + c_0(q) - c_{2n}(q) + d_n
\]
with $(d_n)\in\ell_2^{2\al-1}$, which is to be compared to the
corresponding result of~\cite{KM}, see~\eqref{eq:1.8}. One can also show as
in~\cite{KM} that this estimate is uniform in~$q$ from bounded subsets
of~$W_2^{\al-1}$.
\end{remark}

\begin{remark}\label{rem:Robin}
If $\al=1$, then $\si(x) = h + \int_0^x q$ for some $h\in\bC$, and
formulae~\eqref{eq:1.10} give
\[
    \la_n = \pi n +\frac{\int_0^1 q}{2\pi n} + \frac{\wt\la'_n}{n},
        \qquad
    \wt\mu_n = \pi(n-\tfrac12) + \frac{2h+\int_0^1 q}{2\pi n}
        + \frac{\wt\mu'_n}{n},
\]
for some $\ell_2$-sequences $(\wt\la'_n)$ and $(\wt\mu'_n)$. The choice
$h=0$ corresponds to the genuine Neumann boundary condition $y'(0)=0$
at $x=0$ for the operator $T_{\mathrm N}$, while $h\ne0$ produces the
Robin boundary condition, $y'(0)=hy(0)$, cf.~\eqref{eq:1.5a} and below.
\end{remark}

%%%%%%%%%%%%%%%%%%%%%%%%%%%%%%%%%%%%%%%%%%%%%%%%%%%%%%%%%%%%

\section{Application to inverse spectral problems}\label{sec:appl}

In a selfadjoint regular situation (i.e., for real-valued
integrable~$q$) the classical result of the inverse spectral
theory states that the spectra $\la_n^2$ and $\mu_n^2$ of
Sturm--Liouville operators $T_{\mathrm D}$ and $T_{\mathrm N}$
completely determine the potential. In general the reconstruction
algorithm is quite nontrivial and requires solvability of the
so-called Gelfand--Levitan--Marchenko equation. In this section we
shall show how Theorem~\ref{thm:main0} can be used in the inverse
spectral analysis for some singular (complex-valued) potentials.

Given the eigenvalues $\la_n^2$ and $\mu_n^2$ of $T_{\mathrm D}$ and
$T_{\mathrm N}$ respectively, we denote by $\si^*$ the function
\begin{equation}\label{eq:7.1}
    \si^*(t): = 2\sum_{n=1}^\infty \bigl(\wt\mu_n \sin[(2n-1)\pi t] -
                    \wt\la_n \sin[2\pi nt]\bigr).
\end{equation}
Here, as usual, $\wt\la_n$ and $\wt\mu_n$ are defined
through~\eqref{eq:1.6} and thus the series converges in~$L_2$ as
soon as $q\in W^{\al-1}_2$ with $\al\in[0,1]$. We also observe
that, according to Theorem~\ref{thm:main0}, $\si^* = \si + \wt\si$
for the function $\wt\si\in W^{2\al}_2$ of that theorem.

\begin{proposition}\label{pro:7.1}
Assume that $q\in W^{\al-1}_2$, $\al\in(0,1)$, is such that the
function~$\si^*$ of~\eqref{eq:7.1} belongs to $W^\beta_2$ for some
$\beta\in(0,1]$, $\beta>\al$. Then $q\in W^{\beta-1}_2$.
\end{proposition}

\begin{proof}
We use the so-called bootstrap method. Since $\si^* = \si + \wt\si$ and
$\wt\si\in W^{2\al}_2$ in virtue of Theorem~\ref{thm:main0}, we claim
that $\si$ in fact belongs to $W^{\al'}_2$ with $\al' =
\min\{2\al,\beta\}>\al$. Thus the exponent $\al$ can either be taken
equal to $\beta$ or otherwise doubled. Repeating this procedure
finitely many times, we reach the desired conclusion that $\si\in
W^\beta_2$, i.e., that $q\in W^{\beta-1}_2$.
\end{proof}

Developing the above arguments, we conclude that the function
$\sigma^*$ determines all principal singularities of the
potential~$q$.  The meaning of this claim is that the functions
$\si$ and $\si^*$ share all the singularities typical for
$W^\al_2$. We illustrate this issue by the following example.

Assume that $q$ is such that $\si(x)=\int_0^x q(t)\,dt$ has
bounded variation over $[0,1]$ (i.e., $d\si$ is a finite Borel
measure). Slightly abusing the terminology, we shall say that the
potential $q$ of the corresponding Sturm--Liouville operators
$T_{\si,\mathrm{D}}$ and $T_{\si,\mathrm{N}}$ is a finite Borel
measure; see also~\cite{BR,BEKS,BFT,Zh} for precise definitions.

\begin{proposition}\label{pro:7.2}
Assume that $q$ is a finite Borel measure and $\si^*$ is the
function of~\eqref{eq:7.1} constructed through the corresponding
Dirichlet $\la_n^2$ and Neumann--Dirichlet $\mu_n^2$ eigenvalues.
Then the discrete parts of the measures $q$ and $d\si^*$ coincide.
\end{proposition}

\begin{proof}
Observe that $\si=\int q$, being of bounded variation over
$[0,1]$, belongs to $W_2^{\frac12-\eps}$ for all
$\eps\in(0,\tfrac12)$; henceforth $\wt\si=\si^*-\si$ belongs to
$W^{\frac12+\eps}_2$ for all $\eps\in(0,\tfrac12)$ and thus is
continuous. Therefore $\si^*$ has the same (jump) discontinuities
as $\si$, i.e., the measures~$d\si^*$ and $q=d\si$ have the same
discrete parts.
\end{proof}

As a final remark, we note the following. Roughly speaking,
formula~\eqref{eq:1.10} implies that the Dirichlet spectrum determines
the even part of the potential~$q$, while the Neumann--Dirichlet
spectrum the odd part of~$q$. This explains why the unique
reconstruction of the potential by the Dirichlet spectrum is possible
if the odd part of~$q$ is prespecified, see~\cite{PT}.

%%%%%%%%%%%%%%%%%%%%%%%%%%%%%%%%%%%%%%%%%%%%%%%%%%%%%%%%%%%%
\vspace{.5cm} \noindent\textbf{Acknowledgement. }{The first author
acknowledges the support of the Alexander von Humboldt foundation
and thanks the staff of Institut f\"ur Angewandte Mathematik der
Universit\"at Bonn, in which part of this work was done, for warm
hospitality.}

%%%%%%%%%%%%%%%%%%%%%%%%%%%%%%%%%%%%%%%%%%%%%%%%%%%%%%%%%%%%

\appendix
\section{Interpolation and all that}\label{sec:app.A}

We recall here some facts about interpolation between $W^\al_2$
spaces. For details, we refer the reader to~\cite{BL,KPS,LM}.

 By definition, the space $W^0_2$
coincides with $L_2$ and the norm~$|\cdot|_0$ in $W^0_2$ is just
the $L_2$-norm. The Sobolev space $W^2_2$ consists of all
functions $f$ in $L_2$, whose distributional derivatives $f'$ and
$f''$ also fall into~$L_2$. Being endowed with the norm
\[
    |f|_2 := \bigl(|f|_0^2 + |f'|_0^2 + |f''|_0^2\bigr)^{1/2},
\]
$W_2^2$ becomes the Hilbert space.

Now we interpolate between $W^2_2$ and $W^0_2$ to get the intermediate
spaces $W^s_2$ with norms $|\,\cdot\,|_s$ for $s\in(0,2)$; namely,
$W^{2s}_2 := [W^2_2,W^0_2]_{1-s}$. The norms $|\,\cdot\,|_s$ are
nondecreasing with $s\in[0,2]$, i.e., if $s<t$ and $f\in W^t_2$, then
$|f|_s\le|f|_t$. The space $W^1_2$ consists of $L_2$ functions whose
distributional derivative is again in $L_2$, and
 \(
    \bigl( |f|^2_0 + |f'|^2_0\bigr)^{1/2}
 \)
is the norm on $W^1_2$ that is equivalent to $|\cdot|_1$. Also, the Sobolev
imbedding theorem implies that $W^\al_2$ for $\al>1/2$ is continuously
embedded into $C[0,1]$, the space of continuous functions.

The following result~\cite{BL,KPS,LM} provides a powerful tool in
the study of mappings between the spaces $W_2^s$.

\begin{theorem}[Interpolation Theorem]\label{thm:A.interp}
Assume that $0\le s_0<s_1$, $0\le r_0<r_1$,
 and let $T$ be a linear
operator such that
\[
    |Tf|_{r_0} \le C_0 |f|_{s_0}, \qquad
    |Tg|_{r_1} \le C_1 |g|_{s_1}
\]
for all $f\in W_2^{s_0}$ and all $g\in W_2^{s_1}$. Put $s_t :=
(1-t)s_0 + ts_1$ and $r_t := (1-t)r_0 + tr_1$; then  for every
$t\in[0,1]$ the operator $T$ acts boundedly from $W^{s_t}_2$ to
$W^{r_t}_2$ with norm not exceeding $C_0^{1-t} C_1^t$.
\end{theorem}

We formulate next one result on multilinear interpolation
from~\cite{BL} adapted to our purposes. It concerns analytic
scales of Banach spaces (see also~\cite{KPS}), and we observe that
the scales $\{W_2^{a\theta+b}\}_{\theta\in[0,1]}$ are analytic for
any $a>0$ and $b\ge 0$.

\begin{theorem}\label{thm:polylin-interp}
\def\cJ{{\mathcal J}}
\def\bg{{\mathbf g}}
Assume that $\{E_\theta\}_{\theta\in[0,1]}$ and
$\{G_\theta\}_{\theta\in[0,1]}$ are analytic scales of Banach
spaces, $n\in\mathbb N$, and that $\cJ$ is a multilinear mapping
from $(E_1)^n$ into $G_1$ satisfying the inequalities
\[
    \|\cJ(\bg)\|_{G_0} \le C_0 \prod_{j=1}^n \|g_j\|_{E_0},
    \qquad
    \|\cJ(\bg)\|_{G_1} \le C_1 \prod_{j=1}^n \|g_j\|_{E_1}
\]
for some positive constant $C_0$, $C_1$ and all $\bg:=
(g_1,g_2,\dots,g_n)\in (E_1)^n$. Then $\cJ$ can be uniquely
extended to a multilinear mapping from $(E_\theta)^n$ into
$G_\theta$, $0\le\theta\le1$, with norm not exceeding
$C_0^{1-\theta}C_1^\theta$.
\end{theorem}

The following result is used in various places of the paper.

\begin{lemma}\label{lem:A.2al}
Assume that $f\in W^\al_2$ and $g\in W^\beta_2$ for some
$\al,\beta\in[0,1]$. Then the function \(  h(x) := \int_0^x (fg)
\) belongs to $W^{\al+\beta}_2$ and there exists some constant~$C$
such that
\begin{equation}\label{eq:A.1}
    |h|_{\al+\beta} \le C |f|_\al |g|_\beta.
\end{equation}
\end{lemma}

\begin{proof}
For $\al=\beta=0$ the function $h$ is absolutely continuous and
$|h(x)|\le |f|_0|g_0|$ by the Cauchy--Schwarz inequality, so
that~\eqref{eq:A.1} holds with $C=1$. For $\al=0$ and $\beta=1$,
we find that
\[
    |h'|_0 = |fg|_0 \le (\max_x |g(x)|) |f|_0 \le C_1 |f|_0|g|_1;
\]
henceforth $h\in W^1_2$ and~\eqref{eq:A.1} is satisfied for some
$C=C_2\ge1$. To treat the case $\al=1$, $\beta=0$, just
interchange $f$ and $g$. For $\al=\beta=1$ we find that $h'' = f'g
+ f g' \in L_2$ so that $h\in W^2_2$ with
 \(
    |h|_2 \le C_3 |f|_1 |g|_1.
 \)

Consider now a mapping $M_f:\, L_2 \to L_2$ given by
 \(
        M_fg (x) = \int_0^x (fg)
 \)
with $f\in L_2$ fixed. Then, by the above, $M_f$ acts boundedly in
$L_2$ and $W^1_2$ and
\[
    \|M_f\|_{L_2 \to L_2},\ \|M_f\|_{W^1_2\to W^1_2} \le C_2 |f|_0.
\]
By interpolation, $M_f$ is continuous in $W^\beta_2$ for every
$\beta\in[0,1]$ and its norm $\|M_f\|_{W^\beta_2\to W^\beta_2}$ is
bounded by $C_2|f|_0$; in particular,
\[
    |h|_\beta \le C_2 |f|_0|g|_\beta.
\]
Analogously, for a fixed $f\in W^1_2$ the operator~$M_f$ maps
continuously $W^\beta_2$ into~$W^{1+\beta}_2$ and
\[
    |h|_{1+\beta} \le \max\{C_2,C_3\}\, |f|_1|g|_\beta.
\]
The two above-displayed formulae show that $M_g$ for a fixed $g\in
W^{\beta}_2$ is continuous as a mapping from $L_2$ into
$W^\beta_2$ and from $W^1_2$ into $W^{1+\beta}_2$. Interpolation
now yields its continuity as a mapping from $W^\al_2$ into
$W^{\al+\beta}_2$ for an arbitrary $\al\in[0,1]$ and establishes
inequality~\eqref{eq:A.1}. The lemma is proved.
\end{proof}

The next lemma establishes the required continuity of the mapping
$I_2$, see Section~\ref{sec:smooth}.

\begin{lemma}\label{lem:A.I2}
Assume that $f\in W^\al_2$ and $g\in W^\beta_2$ for some
$\al,\beta\in[0,1]$. Then the function
\[
  h(x) := \int_0^{1-x} f(x+t)g(t)\,dt
\]
belongs to $W^{\al+\beta}_2$ and there exists some constant~$C$
such that
\[
    |h|_{\al+\beta} \le C |f|_\al |g|_\beta.
\]
\end{lemma}

\begin{proof}
The proof is completely analogous to the one of
Lemma~\ref{lem:A.2al}: one establishes first the statement for
$\al$ and $\beta$ equal to $0$ or $1$, and then interpolate. The
only remark is that for $\beta=1$ one should use a representation
 \(
    h(x) = \int_x^1 f(s) g(s-x)\,ds
 \)
to be able to differentiate $h$.
\end{proof}

\begin{lemma}\label{lem:A.5}
{\rm(a)} Assume that $f\in W^\al_2$, $\al\in[0,1]$; then
 $s_{2n}(f) \in\ell_\infty^\al$ and  $c_{2n}(f) \in \ell_2^\al$.
{\rm(b)} Assume that $g\in W^\beta_2$, $\beta\in (1,2]$, and
$g(0)=g(1)$; then $(s_{2n}(g))\in\ell_2^\beta$.
\end{lemma}

\begin{proof}
Part (a) follows by interpolation. Indeed, the mapping
\[
      S: f\mapsto (s_{2n}(f))_{n\in\mathbb N}
\]
is bounded from $L_2$ into $\ell_\infty$ and from $W_2^1$ into
$\ell_\infty^1$. Since the spaces $\ell_\infty^s$ form an analytic
Banach scale~\cite{KPS}, Interpolation
Theorem~\ref{thm:polylin-interp} yields the result about
$s_{2n}(f)$. Similar arguments justify the statement about
$c_{2n}(f)$.

Part (b) requires only a slight modification. For an arbitrary $s\in [1,2]$
we put
\[
     \wt W_2^s :=\{h\in W_2^s \mid h(0)=h(1) \}.
\]
Since the family $\{W_2^{\theta +1}\}_{\theta\in[0,1]}$ forms a
Hilbert scale, by virtue of the general interpolation result
from~\cite[Ch.1, 13.4]{LM}, the family
 $\{\widetilde W_2^{\theta +1}\}_{\theta\in[0,1]}$
also is a Hilbert scale. Simple integration by parts shows that
the above operator~$S$ maps continuously $\wt W_2^1$ into
$\ell_2^1$ and $\wt W_2^2$ into $\ell_2^2$. By Interpolation
Theorem~\ref{thm:polylin-interp} the operator $S$ maps
continuously $\wt W_2^{\theta +1}$ into
 $\ell_2^{\theta+1}$ for all $\theta\in[0,1]$, and the proof is complete.
\end{proof}

\begin{remark}\label{rem:A.sin-coef}
Part (b) of the lemma also holds true for $\beta\in(\tfrac12,1)$
(see, e.g.,~\cite{KM}).
\end{remark}

For $L_2$-functions $f$ and $g$ we denote by $f\ast g$ their convolution,
i.e.,
\[
    (f \ast g)(x) := \int_0^x f(x-t)g(t)\,dt.
\]
Recall also that $R$ is the reflection operator about $x=\tfrac12$:
$(Rf)(x) = f(1-x)$.

\begin{lemma}\label{lem:A.3}
Suppose that $f,g\in L_2$; then we have
\[
    c_n(f) c_n(g) = c_n (h_1), \qquad
    s_n(f) s_n(g) = c_n (h_2), \qquad
    s_n(f) c_n(g) = s_n (h_3),
\]
where the functions $h_j$, $j=1,2,3$, are given by
\begin{equation}\label{eq:A.13}
\begin{aligned}
    h_1 &= \frac12\bigl\{R(Rf\ast g + f\ast Rg) +
        f\ast g + Rf\ast Rg\bigr\},\\
    h_2 &= \frac12\bigl\{R(Rf\ast g + f\ast Rg) -
        f\ast g - Rf\ast Rg\bigr\},\\
    h_3 &= \frac12\bigl\{R(Rf\ast g - f\ast Rg) +
        f\ast g - Rf\ast Rg\bigr\}.
\end{aligned}
\end{equation}
\end{lemma}

\begin{proof}
We shall prove only the first equality, since the other ones are treated
analogously. We have
\[
    2 c_n(f)c_n(g) = \int_0^1 \int_0^1 f(x) g(t)
            [\cos\pi n(x-t) + \cos\pi n(x+t)]\,dx dt,
\]
and simple calculations lead to
\begin{align*}
    \int_0^1\int_0^1 &f(x) g(t)\cos\pi n(x-t)\, dx dt \\
        &= \int_0^1
        \Bigl( \int_0^{1-s} f(s+t) g(t)\,dt
                + \int_0^{1-s} f(t) g(s+t)\,dt
        \Bigr)\cos\pi ns\,ds, \\
    \int_0^1\int_0^1 &f(x) g(t)\cos\pi n(x+t)\, dx dt \\
        &= \int_0^1
        \Bigl( \int_0^{s} f(s-t) g(t)\,dt
                + \int_0^{s} f(1-t) g(1-s+t)\,dt
        \Bigr)\cos\pi ns\,ds.
\end{align*}
Taking into account the relations
\[
    \int_0^{1-s} f(s+t) g(t)\,dt = R(Rf\ast g)(s), \qquad
    \int_0^{s} f(1-t) g(1-s+t)\,dt = Rf \ast Rg,
\]
we get $c_n(f) c_n(g) = c_n(h_1)$ with $h_1$ as stated. The lemma is
proved.
\end{proof}

\begin{proof}[Proof of Lemma~\ref{lem:5.2}]
By Lemma~\ref{lem:A.3} we have that $\bc(f)\bc(g)=\bc(h_1)$,
$\bs(f)\bs(g)=\bc(h_2)$, and $\bs(f)\bc(g)=\bs(h_3)$ for the
functions $h_1$, $h_2$, and $h_3$ given by~\eqref{eq:A.13}.

It is easily verified that, with the operator $I_2$ of
Section~\ref{sec:smooth}, we have
\[
        I_2(f,g) = R(Rf\ast g).
\]
Recalling that $R$ is unitary in every $W^\al_2$, we conclude by
Lemma~\ref{lem:A.I2} that the functions $h_j$ of~\eqref{eq:A.13}
are in $W^{\al+\beta}_2$ as soon as $f\in W^\al_2$ and $g\in
W^\beta_2$ for some $\al,\beta\in[0,1]$ and that, moreover, with
some $\rho>0$ the inequality
\[
    |h_j|_{\al+\beta} \le \rho |f|_\al|g|_\beta
\]
holds. Recalling the definition of the norm in $\mS_\al$ and
$\mC_\al$, we get the required estimates~\eqref{eq.51}.
\end{proof}

%%%%%%%%%%%%%%%%%%%%%%%%%%%%%%%%%%%%%%%%%%%%%%%%%%%%%%%%%%%%

\end{document}